# High Dimensional Three-Periods Locally Ideal MIP Formulations for the UC Problem

Linfeng Yang, *Member*, *IEEE*, Wei Li, Yan Xu, *Senior Member*, *IEEE*, Cuo Zhang, Beihua Fang

*Abstract*—The thermal unit commitment (UC) problem often can be formulated as a mixed integer quadratic programming (MIQP), which is difficult to solve efficiently, especially for large-scale instances. The tighter characteristic reduces the search space, therefore, as a natural consequence, significantly reduces the computational burden. In the literature, many tightened formulations for single units with parts of constraints were reported without presenting how they were derived. In this paper, a systematic approach is developed to formulate the tight formulations. The idea is using more new variables in high dimension space to capture all the states for single units within three periods, and then, using these state variables systematic derive three-periods locally ideal expressions for a subset of the constraints in UC. Meanwhile, the linear dependence relations of those new state variables are leveraged to keep the compactness of the obtained formulations. Based on this approach, we propose two tighter models, namely 3P-HD and 3P-HD-Pr. The proposed models and other four state-of-the-art models were tested on 51 instances, including 42 realistic instances and 9 8-unit-based instances, over a scheduling period of 24 h for systems ranging from 10 to 1080 generating units. The simulation results show that our proposed MIQP UC formulations are the tightest and can be solved most efficiently. After using piecewise technique to approximate the quadratic operational cost function, the six UC MIQP formulations can be approximated by six corresponding mixed-integer linear programming (MILP) formulations. Our experiments show that the proposed 3P-HD and 3P-HD-Pr MILP formulations also perform the best in terms of tightness and solution times.

NOMENCLATURE

Indices:
$i$      Index for unit.
$t$      Index for time period.
Operator:
$[\cdot]^+$      $\max(0,\cdot)$
$\gtrsim (\eqsim, \lesssim) \; F_1 \gtrsim (\eqsim, \lesssim) F_2$ means that MIP formulation (or constraints) $F_1$ is tighter than (equivalent to, looser than) $F_2$ in tightness.
Constants:
$N$      Total number of units.
$T$      Total number of time periods.
$\mathbb{N}$      Set of indexes of units.
$\alpha_i, \beta_i, \gamma_i$      Coefficients of the quadratic production cost function of unit $i$.
$C_{\text{hot},i}$      Hot startup cost of unit $i$.
$C_{\text{cold},i}$      Cold startup cost of unit $i$.
$\underline{T}_{\text{on},i}$      Minimum up time of unit $i$.
$\underline{T}_{\text{off},i}$      Minimum down time of unit $i$.
$T_{\text{cold},i}$      Cold startup time of unit $i$.
$\overline{P}_i$      Maximum power output of unit $i$.
$\underline{P}_i$      Minimum power output of unit $i$.
$P_{\text{D},t}$      System load demand in period $t$.
$R_t$      Spinning reserve requirement in period $t$.
$P_{\text{up},i}$      Ramp up limit of unit $i$.
$P_{\text{down},i}$      Ramp down limit of unit $i$.
$P_{\text{start},i}$      Startup ramp limit of unit $i$.
$P_{\text{shut},i}$      Shutdown ramp limit of unit $i$.
$u_{i,0}$      Initial commitment state of unit $i$ (1 if it is online, 0 otherwise).
$T_{i,0}$      Number of periods unit $i$ has been online ($+$) or offline ($-$) prior to the first period of the time span (end of period 0).
$U_i$      $\left[\min[T, u_{i,0}(\underline{T}_{\text{on},i} - T_{i,0})]\right]^+$
$L_i$      $\left[\min[T, (1-u_{i,0})(\underline{T}_{\text{off},i} + T_{i,0})]\right]^+$
Variables:
$u_{i,t}$      Schedule of unit $i$ in period $t$, binary variable that is equal to 1 if unit $i$ is online in period $t$ and 0 otherwise.
$s_{i,t}$      Startup status of unit $i$ in period $t$.
$d_{i,t}$      Shutdown status of unit $i$ in period $t$.
$P_{i,t}$      Power output of unit $i$ in period $t$.
$S_{i,t}$      Startup cost of unit $i$ in period $t$.

This work was supported by the Natural Science Foundation of China (51767003), the Guangxi Natural Science Foundation (2017GXNSFBA198238). (Corresponding author: Yan Xu.)
L.F. Yang, W. Li and B.H. Fang are with the School of Computer Electronics and In-formation, Guangxi University, Nanning 530004, China. And L.F. Yang is also with the Guangxi Key Laboratory of Multimedia Communication and Network Technology, Guangxi University. W. Li is also with the School of Electrical Engineering, Guangxi University. (e-mail: ylf@gxu.edu.cn; 2863053347@qq.com; 13706108501@163.com).
Y. Xu is with School of Electric and Electronic Engineering, Nanyang Technological University, Singapore. (e-mail: xuyan@ntu.edu.sg)
C. Zhang is with School of Electric Engineering and Telecommunications, University of New South Wales, Sydney, Australia.

## I. INTRODUCTION

THE unit commitment (UC) continues to attract significant attention from both industry and academia because it is an extremely important optimization problem for both daily operation scheduling and planning studies from short term to long term. In general, the UC problem is formulated as a mixed integer nonlinear programming (MINLP) problem [1]-[3] to determine the operation schedule of the generating units at each time period with varying loads under different operating constraints and environments. However, the large-scale nature and nonconvexity of the problem make it challenging to solve. Consequently, developing solution methods that can achieve high-quality solutions in a short amount of time has been the focus of significant research over the last several decades. Many reported methods, including artificial intelligence (AI) algorithms [4]-[6] and numerical optimization techniques [7]-[26], have been developed for solving the UC problem. The AI methods include evolutionary algorithms [4], particle swarm optimization [5], simulated annealing [6] and so on. The numerical optimization methods include priority list (PL) algorithms [7], outer approximation (OA) [8], [9], Lagrangian relaxation (LR) [10], and MIP methods [11]-[26]. AI methods can produce fair solutions within a reasonable computation time. However, the quality of the solutions is difficult to guarantee [1]. Among the numerical optimization methods, PL algorithms are the easiest to implement with a fast convergence rate. However, such algorithms usually suffer from being highly heuristic in nature and yield relatively poor-quality solutions [7]. OA approaches can theoretically produce exact (or nearly exact) solutions. However, these algorithms have impractical computation time requirements for large-scale systems [8], [9]. Historically, LR has been the method of choice for the UC scheduling software that is used in the power industry for large-scale instances when a very fast computation time is a priority [10]. However, in general, the solution obtained via LR by solving the dual problem of the UC problem is not feasible. Therefore, heuristic approaches are needed to assist in searching for feasible solutions. This sub-optimality is LR's main shortcoming. For a detailed review, the reader is referred to [1].

With the significant progresses in the theory of mixed integer programming (MIP) and improvements in the efficiency of general-purpose branch-and-cut solvers in recent years [27]-[31], solving UC problems by using MIP method is becoming increasingly popular. In 2005, the world's largest competitive wholesale market PJM and other power systems operators in the United States began to change from LR to MIPs to tackle its UC-based scheduling problems. The savings associated with this transition are estimated at $5 billion USD annually in the United States alone [13].

It is well-known that the choice of MIP formulation for any given optimization problem can have a significant impact on practical computational difficulty [29], [31], and UC is no exception. Consequently, there has been significant research focus over the past twelve years in developing improved UC MIP formulations, i.e., focusing on "tighter" or "compacter" formulations. [11] proposed a way to approximating the nonlinear objective function based on perspective-cut. [2] presents a more compact formulation by reducing two sets of binary variables from the three binary variables formulation of [14]. Moving forward, tight and compact UC formulation is provided in [15], but this paper mainly focuses on the reformulation of the startup costs of units and does not consider the nonlinear production cost. [12] proposed a novel two-binary-variable MIQP formulation for the UC problem. Compared with one-binary-variable and three-binary-variable formulations, the 2-bin formulation is more compact, and it is tighter in terms of the quadratic cost function. The tighter UC models are the more studied. Most of the UC literature involves finding a *locally ideal* or locally tighter formulation for a subset of the constraints in UC models, i.e., for one unit, just considering one or several of the types of constraints in two or several periods and deriving a result for that or those constraint type(s). For instance, [16] have identified the convex hull for a minimum up/downtime polytope, [17] provided the convex hull of generation limits and this minimum up/downtime polytope, [18] provided the convex hull of a two-period ramping polytope and exponential classes of multi-period variable upper bound and multi-period ramping facets for some giving conditions, [19] provided strengthened inequalities of a ramping polytope. We will give detailed explanations for these tight constraints in Section II of this paper and the reader is referred to [20] for a more detailed review.

Nevertheless, it is well known that in general, MIPs are NP-hard, and solving UC problems of realistic size, involving thousands of generators, over several time periods, remains challenging [21]-[22]. Moreover, improving MIP formulations can dramatically reduce its computational burden and so allow the implementation of more advanced and computationally demanding problems [23], such as valve-point effect [24], stochastic formulation [25], or transmission switching [26].

In this paper, we study MIP formulations of the UC problem with components that are of critical interest. The contribution of our study is threefold:

1) First, by introducing more state variables, we present the basic procedure for deriving *locally ideal* formulations for generations limits and production ramping limits in high dimension space while three time periods are considered. And after projecting out or conditional dropping the introduced auxiliary state variables, our high dimension constraints can easily imply the other state-of-the-art constraints, most of these formulations are constructed rely on experiences in original literatures. In addition, the proposed high dimension constraints also imply some new local ideal constraints in classical UC variable space.

2) Second, a tight and compact high dimension UC model is constructed and a new UC model in the classical variable space is deduced by eliminating the auxiliary state variables.

3) Finally, we perform a computational study to demonstrate the behaviors of the two news and several relatively recent and well-regarded benchmark formulations.

The remaining parts of this paper are organized as follows. In Section II, the state-of-the-art MIQP formulations of the UC problem are reviewed. In section III, we present the basic procedure for deriving our locally ideal formulations and show their relations to the other classical constraints. The computational results are reported and analyzed in Section IV to verify the effectiveness of the proposed formulations and models. Finally, we conclude the paper in Section V. Meanwhile, some proofs of theorems and experimental results are shown in Appendix.

## II. UC Problem and its State-of-the-art Formulations

In general, the UC problem can be formulated as follows [20]
$$\min F_C = \sum_{i=1}^{N} \sum_{t=1}^{T} [f_i(P_{i,t}) + S_{i,t}]$$
$$\text{s.t.} \begin{cases} \sum_{i=1}^{N} A_i(\boldsymbol{P}_i, \boldsymbol{y}_i) + N(\boldsymbol{z}) = L \\ (\boldsymbol{P}_i, \boldsymbol{y}_i, \boldsymbol{S}_i) \in \Pi_i. \end{cases} \quad (1)$$

where $F_C$ is the total operation cost, the production cost is $f_i(P_{i,t}) = \alpha_i u_{i,t} + \beta_i P_{i,t} + \gamma_i (P_{i,t})^2$, and $S_{i,t}$ is the startup cost. In what follows, a vector of variables of the same type (say $P_{i,t}$, $\forall t \in \{1, ..., T\}$) will be typed in bold ($\boldsymbol{P}_i$). $\boldsymbol{y}_i$ is the vector including all status variables (such as $\boldsymbol{u}_i$, $\boldsymbol{s}_i$, $\boldsymbol{d}_i$, etc.) for generator $i$. The matrix $A_i(\boldsymbol{P}_i, \boldsymbol{y}_i)$ determines how the generator interacts with the system requirements, which are written in matrix form as the first constraint in (1). Here the variable $\boldsymbol{z}$ is other potential decision variables involving the operation of the system. Finally, the last constraint in (1) defines the constraints for each generator's schedule, often includes minimum and maximum generation levels when on, production ramping limits, minimum up-times and down-times, time-dependent start-up costs. Although the set $\Pi_i$ only includes linear constraints, it is non-convex because of the binary variables $\boldsymbol{y}_i$.

Now we consider recent UC formulations from the literatures. When describing the formulation for a single generator, we sometimes drop the subscript $i$ for clarity and ease of presentation. Most of the constraints described here could be applied to every generator $i \in \{1, ..., N\}$ and every period $t \in \{1, ..., T\}$ if such indexes make any sense, and we will explicitly point out when this is not the case.

### A. System Constraints
*1) Power balance constraint:*
$$\sum_{i=1}^{N} P_{i,t} - P_{D,t} = 0. \quad (2)$$
In order to maintain system security, the total power output of generating units must satisfies total load demand.

*2) System spinning reserve requirement:*
$$\sum_{i=1}^{N} u_{i,t} \overline{P}_i \geq P_{D,t} + R_t. \quad (3)$$
Spinning reserve is an important resource, which is used by system operators to maintain system security and improve operational reliability.

### B. Constraints and variables describing $\Pi_i$
*1) State variables and logical constraints*
The prototypical MIP formulation of UC [14] use three binary variables to represent the statuses of unit, i.e., on/off ($u_t$), startup ($s_t$), shutdown ($d_t$). [14] also firstly formulated the logical constraints to relate the three binary variables:
$$s_t - d_t = u_t - u_{t-1}. \quad (4)$$
[2] presents the first 1-bin UC formulation (with only $u_t$) by reducing two sets of binary variables from the 3-bin model. However, the assumed benefit of using considerably smaller number of binary variables does not necessarily lead to superior computational performance, as observed by [19]. This is due to improved 3-bin formulations, especially in tightness, and more robust MIP solvers.

In addition, one may use equation (4) to project out either the $s_t$ or $d_t$ variables while not losing strength. If the compactness of the model has been improved in this projecting process, MIP solvers will obtain better performance. For example, [12] completely projects out the shutdown variables $d_t$.

[32] suggests replacing the variable $u_t$ with a state-transition variable $o_t$ ($\tilde{u}_t$ was used in [32]) which encodes if the generator remains operational at time $t$, i.e., $u_{t-1} = u_t = 1$. They also give the mathematical relationship,
$$u_t = o_t + s_t, \quad (5)$$
which allows for the transformation of classical 3-bin model to the state-transition formulation.

*2) minimum up/down time constraints*
In [16], Rajan and Takriti use $s_t$ and $d_t$ to formulate the following minimum up/down time constraints
$$\sum_{\varpi=[t-\underline{T}_{\text{on},i}]^+ +1}^{t} s_\varpi \leq u_t, \ t \in [U_i + 1, ..., T], \quad (6)$$
$$\sum_{\varpi=[t-\underline{T}_{\text{off},i}]^+ +1}^{t} d_\varpi \leq 1 - u_t, \ t \in [L_i + 1, ..., T]. \quad (7)$$
Indeed, [16] shows that constraints (6)-(7) define facets of the minimum up/down time polytope, which together with (4) and variable bounds are an ideal formulation for up-time and down-time.

*3) Unit generation limits and upper bounds:*
The simplest generation limit is given by [2]
$$u_{i,t} \underline{P}_i \leq P_{i,t}, \quad (8)$$
$$P_{i,t} \leq u_{i,t} \overline{P}_i. \quad (9)$$
Additionally, to ensure total capacity is not exceeded as the generator $i$ is shutting down, [33] define the constraint
$$P_{i,t} \leq u_{i,t} \overline{P}_i - d_{i,t+1}(\overline{P}_i - P_{\text{shut},i}), \ t \in \{1, ..., T-1\}. \quad (10)$$
[15] proposes using the start-up and shutdown ramping limits to tighten the variable upper bounds for unit $i \in \mathcal{I}^{\geq 2}$:
$$P_{i,t} \leq u_{i,t} \overline{P}_i - s_{i,t}(\overline{P}_i - P_{\text{start},i}) - d_{i,t+1}(\overline{P}_i - P_{\text{shut},i}). \quad (11)$$
And for $i \in \mathcal{I}^1$:
$$P_{i,t} \leq u_{i,t} \overline{P}_i - s_{i,t}(\overline{P}_i - P_{\text{start},i}), \quad (12)$$
$$P_{i,t} \leq u_{i,t} \overline{P}_i - d_{i,t+1}(\overline{P}_i - P_{\text{shut},i}). \quad (13)$$
where $\mathcal{I}^{\geq 2} \coloneqq \{i | \underline{T}_{\text{on},i} \geq 2\}$ and $\mathcal{I}^1 \coloneqq \{i | \underline{T}_{\text{on},i} = 1\}$.

[18] showed that (12) and (13) are facets of the two-period ramp-up and ramp-down polytopes, respectively.

For $i \in \mathcal{I}^1$, [17] proposes the following upper bounds constraints:
$$P_{i,t} \leq u_{i,t} \overline{P}_i - s_{i,t}(\overline{P}_i - P_{\text{start},i}) - d_{i,t+1}[P_{\text{start},i} - P_{\text{shut},i}]^+. \quad (14)$$
$$P_{i,t} \leq u_{i,t} \overline{P}_i - d_{i,t+1}(\overline{P}_i - P_{\text{shut},i}) - s_{i,t}[P_{\text{shut},i} - P_{\text{start},i}]^+. \quad (15)$$

And [17] shows that (14)-(15) are strict tighter than (12)-(13) when $P_{\text{shut},i} \neq P_{\text{start},i}$. Furthermore, [17] shows that (11) and (14)-(15) all are the cornerstone constraints which form the convex hull description for the following basic operating constraints of a single unit in multiple periods: 1) generation limits, 2) startup and shutdown capabilities, and 3) minimum up/down times.

*4) Ramping constraints*

The popular ramping constraints are given by [33] as

$$P_t - P_{t-1} \leq u_{t-1}P_{\text{up}} + s_t P_{\text{start}} \tag{16}$$
$$P_{t-1} - P_t \leq u_t P_{\text{down}} + d_t P_{\text{shut}}. \tag{17}$$

[18] provided the strengthened ramp up/down constraints

$$P_t - P_{t-1} \leq u_t(P_{\text{up}} + \underline{P}) - u_{t-1}\underline{P} + s_t(P_{\text{start}} - P_{\text{up}} - \underline{P}), \tag{18}$$

$$P_{t-1} - P_t \leq u_{t-1}(P_{\text{down}} + \underline{P}) - u_t\underline{P} + d_t(P_{\text{shut}} - P_{\text{down}} - \underline{P}), \tag{19}$$

which were proved to be facet-defining for the two-period ramp rate polytopes [18].

[19] introduced strengthened inequalities for ramping, under certain assumptions on the generator. A ramp-up inequality is proposed for unit $i \in \mathcal{J}^{\geq 2} \cap \mathcal{L}$

$$P_t - P_{t-1} \leq u_t P_{\text{up}} - d_t \underline{P} - d_{t+1}(P_{\text{up}} - P_{\text{shut}} + \underline{P}) + s_t(P_{\text{start}} - P_{\text{up}}), \tag{20}$$

where $\mathcal{L} = \{i | P_{\text{up},i} > P_{\text{shut},i} - \underline{P}_i\}$. And the other ramping constraint, bounded on three periods, is for unit $i \in \mathcal{J}^{\geq 2} \cap \underline{\mathcal{J}}^{\geq 2} \cap \mathcal{L}$ with

$$P_{t+1} - P_{t-1} \leq 2u_{t+1}P_{\text{up}} - d_t\underline{P} - d_{t+1}\underline{P} + s_t(P_{\text{start}} - P_{\text{up}}) + s_{t+1}(P_{\text{start}} - 2P_{\text{up}}), \tag{21}$$

where $\underline{\mathcal{J}}^{\geq 2} \coloneqq \{i | \underline{T}_{\text{off},i} \geq 2\}$.

And a ramp-down inequality for unit $i \in \mathcal{J}^{\geq 2} \cap \underline{\mathcal{L}}$ is

$$P_{t-1} - P_t \leq u_t P_{\text{down}} + d_t P_{\text{shut}} - s_{t-1}(P_{\text{down}} - P_{\text{start}} + \underline{P}) - s_t(P_{\text{down}} + \underline{P}), \tag{22}$$

where $\underline{\mathcal{L}} = \{i | P_{\text{down},i} > P_{\text{start},i} - \underline{P}_i\}$.

[19] shows that (22) is facet of the subspace formed by projecting $\Pi_i$ onto the set of variables $\mathcal{S}_{(i,t-1,t)} = \{P_{i,\tau}, u_{i,\tau}, s_{i,\tau}, d_{i,\tau}: \tau = t-1, t\}$, (20) and (21) are facets of projected subspace on $\mathcal{S}_{(i,t-1,t,t+1)} = \{P_{i,\tau}, u_{i,\tau}, s_{i,\tau}, d_{i,\tau}: \tau = t-1, t, t+1\}$.

Here, we note that, the other generation upper bounds constraints and ramping constraints proposed in [19] and [20] involving unit states more than three periods have not been covered here, mainly because these constraints have not been improved in our work and are equally effective for our model [34]. [18] and [35] also give several exponential classes of variable upper bound inequalities and two-period ramping inequalities, however, such inequalities require separation and hence will not be covered here.

*5) startup cost*

The startup cost $S_{i,t}$ can be represented as an MILP formulation [36]:

$$S_{i,t} \geq C_{\text{hot},i} s_{i,t}, \tag{23}$$

$$S_{i,t} \geq C_{\text{cold},i} \left[ s_{i,t} - \sum_{\tau=\max(t-\underline{T}_{\text{off},i}-T_{\text{cold},i},1)}^{t-1} d_{i,\tau} - f_{\text{init},i,t} \right], \tag{24}$$

where $f_{\text{init},i,t} = 1$ when $t - \underline{T}_{\text{off},i} - T_{\text{cold},i} \leq 0$ and $[-T_{i,0}]^+ < |t - \underline{T}_{\text{off},i} - T_{\text{cold},i} - 1| + 1$, $f_{\text{init},i,t} = 0$ otherwise.

Because that $C_{\text{cold},i} \geq C_{\text{hot},i}$ is common, then after introducing $\tilde{S}_{i,t}$ representing the part of startup cost exceeding $C_{\text{hot},i}$, [32], [12] reformulate the objective function as

$$\min F_{\text{C}} = \sum_{i=1}^{N} \sum_{t=1}^{T} [u_{i,t} f_i(P_{i,t}) + C_{\text{hot},i} s_{i,t} + \tilde{S}_{i,t}]. \tag{25}$$

In the same time, (23) can be dropped and (24) should be reformulated as

$$\tilde{S}_{i,t} \geq (C_{\text{cold},i} - C_{\text{hot},i}) \left[ s_{i,t} - \sum_{\tau=\max(t-\underline{T}_{\text{off},i}-T_{\text{cold},i},1)}^{t-1} d_{i,\tau} - f_{\text{init},i,t} \right]. \tag{26}$$

[32] and [12] point out that (26) is compacter than (23)~(24), however, neither of them notice the new formulation for startup cost is "*tighter*" than (23)~(24) because $C_{\text{hot},i} s_{i,t} + \tilde{S}_{i,t} \geq \max\left\{ C_{\text{hot},i} s_{i,t}, C_{\text{cold},i} \left[ s_{i,t} - \sum_{\tau=\max(t-\underline{T}_{\text{off},i}-T_{\text{cold},i},1)}^{t-1} d_{i,\tau} - f_{\text{init},i,t} \right] \right\}$. The detailed proof is given as follows.

**Proof**: We use $\xi_{it}$ to denote $\sum_{\tau=\max(t-\underline{T}_{\text{off},i}-T_{\text{cold},i},1)}^{t-1} d_{i,\tau} + f_{\text{init},i,t}$. It is obvious that $0 \leq \xi_{it} \in \mathbb{Z}$. Then, we use $\Psi_1$ to denote the set $\left\{ S_{i,t} \middle| \begin{array}{c} S_{i,t} \geq C_{\text{hot},i} s_{i,t} \\ S_{i,t} \geq C_{\text{cold},i}[s_{i,t} - \xi_{it}] \\ S_{i,t} \geq 0 \\ s_{i,t} \in \{0,1\} \\ 0 \leq \xi_{it} \in \mathbb{Z} \end{array} \right\}$, $\Psi_2$ to denote $\left\{ S_{i,t} \middle| \begin{array}{c} S_{i,t} = C_{\text{hot},i} s_{i,t} + \tilde{S}_{i,t} \\ \tilde{S}_{i,t} \geq (C_{\text{cold},i} - C_{\text{hot},i})[s_{i,t} - \xi_{it}] \\ \tilde{S}_{i,t} \geq 0 \\ s_{i,t} \in \{0,1\} \\ 0 \leq \xi_{it} \in \mathbb{Z} \end{array} \right\}$. It is not difficult to verify that $\Psi_1 = \Psi_2$.

We use $R(\cdot)$ to denote the continuous relaxation of set ".", now, we will show that $R(\Psi_1) \supseteq R(\Psi_2)$ and for $s_{i,t} > \xi_{i,t} > 0$, $R(\Psi_1) \supset R(\Psi_2)$.

Using (26), we have

$$C_{\text{hot},i} s_{i,t} + \tilde{S}_{i,t} \geq C_{\text{cold},i}(s_{i,t} - \xi_{i,t}) + C_{\text{hot},i} \xi. \tag{27}$$

Since $C_{\text{hot},i} \xi \geq 0$, then we have

$$C_{\text{hot},i} s_{i,t} + \tilde{S}_{i,t} \geq C_{\text{cold},i}(s_{i,t} - \xi_{i,t}). \tag{28}$$

In addition, since $\tilde{S}_{i,t} \geq 0$, we immediately have

$$C_{\text{hot},i} s_{i,t} + \tilde{S}_{i,t} \geq C_{\text{hot},i} s_{i,t} \tag{29}$$

According to (28)-(29), we have $R(\Psi_1) \supseteq R(\Psi_2)$.

When $\xi_{i,t} > 0$, considering (27), we have

$$C_{\text{hot},i} s_{i,t} + \tilde{S}_{i,t} > C_{\text{cold},i}(s_{i,t} - \xi_{i,t}) \tag{30}$$

When $s_{i,t} > \xi_{i,t}$, considering (26), we have $\tilde{S}_{i,t} > 0$, furthermore, we have

$$C_{\text{hot},i} s_{i,t} + \tilde{S}_{i,t} > C_{\text{hot},i} s_{i,t} \tag{31}$$

So, when $s_{i,t} > \xi_{i,t} > 0$, according to (30)-(31), we have $C_{\text{hot},i} s_{i,t} + \tilde{S}_{i,t} > \max\{C_{\text{hot},i} s_{i,t}, C_{\text{cold},i}[s_{i,t} - \xi_{i,t}]\}$. Then, $R(\Psi_1) \supset R(\Psi_2)$. ∎

*6) Initial status of units*

The constraints to enforce the initial uptime and downtime based on the generator's history is

$$u_{i,t} = u_{i,0}, \quad t \in [1, \ldots, U_i + L_i]. \tag{32}$$

*C. State-of-the-art UC formulations*

The tighter characteristic reduces the search space and the more compact characteristic increases the searching speed with which solvers explore that reduced space. However, in the procedure of modeling, tightness and compactness always cannot be satisfied simultaneously. This dilemma can be shown in Section II.B.3) and II.B.4) indeed, a tight constraint always sacrifices some compactness.

When compactness became a priority, the state-of-the-art 3-bin UC MIQP formulation with the most compact constraints for unit generation limits and ramping limits within two periods $(t-1, t)$ [2][33], denoted as 2-period-compact model (2P-Co), is

$$\min F_C = \sum_{i=1}^{N}\sum_{t=1}^{T}\left[\alpha_i u_{i,t} + \beta_i P_{i,t} + \gamma_i (P_{i,t})^2 + S_{i,t}\right]$$

$$\text{s.t.}\begin{cases} (23)(24)(4)(32)(6)(7) \\ (8)(9)(16)(17) \\ (2)(3) \\ (u_{i,t}, s_{i,t}, d_{i,t}) \in \{0,1\}^3, (P_{i,t}, S_{i,t}) \in \mathcal{R}_+^2. \end{cases} \quad (33)$$

Considering compactness and tightness simultaneously, the 3-bin UC formulation with fact-define constraints for unit generation limits and ramping limits within two periods [18], denoted as 2-period-tight model (2P-Ti), is

$$\min F_C = \sum_{i=1}^{N}\sum_{t=1}^{T}\left[\alpha_i u_{i,t} + \beta_i P_{i,t} + \gamma_i (P_{i,t})^2 + S_{i,t}\right]$$

$$\text{s.t.}\begin{cases} (23)(24)(4)(32)(6)(7) \\ (9) \text{ for } t \in \{1,T\}; (8) \\ (9) \text{ for } i \in \mathcal{I}^{\geq 2}; (12)(13) \text{ for } i \in \mathcal{I}^1 \\ (18)(19) \\ (2)(3) \\ (u_{i,t}, s_{i,t}, d_{i,t}) \in \{0,1\}^3, (P_{i,t}, S_{i,t}) \in \mathcal{R}_+^2. \end{cases} \quad (34)$$

When considering three periods, the tightest 3-bin UC formulation with tightest constraints for unit generation limits and ramping limits within three periods [15][19], denoted as 3-period-tight model (3P-Ti), is

$$\min F_C = \sum_{i=1}^{N}\sum_{t=1}^{T}\left[\alpha_i u_{i,t} + \beta_i P_{i,t} + \gamma_i (P_{i,t})^2 + S_{i,t}\right]$$

$$\text{s.t.}\begin{cases} (23)(24)(4)(32)(6)(7) \\ (9) \text{ for } t \in \{1,T\}; (8) \\ (11) \text{ for } i \in \mathcal{I}^{\geq 2}; (14)(15) \text{ for } i \in \mathcal{I}^1, \\ (20) \text{ for } i \in \mathcal{I}^{\geq 2} \cap \mathcal{L}; (22) \text{ for } i \in \mathcal{I}^{\geq 2} \cap \underline{\mathcal{L}} \\ (18) \text{ for } i \in (\mathbb{N} - (\mathcal{I}^{\geq 2} \cap \mathcal{L})) \cup ((\mathcal{I}^{\geq 2} \cap \mathcal{L}) \wedge (t = T)) \\ (19) \text{ for } i \in (\mathbb{N} - (\mathcal{I}^{\geq 2} \cap \underline{\mathcal{L}})) \cup ((\mathcal{I}^{\geq 2} \cap \underline{\mathcal{L}}) \wedge (t = 2)) \\ (21) \text{ for } i \in \mathcal{I}^{\geq 2} \cap \underline{\mathcal{I}}^{\geq 2} \cap \mathcal{L} \\ (2)(3) \\ (u_{i,t}, s_{i,t}, d_{i,t}) \in \{0,1\}^3, (P_{i,t}, S_{i,t}) \in \mathcal{R}_+^2. \end{cases}$$
$$(35)$$

Taking (5) under consideration, 3P-Ti model can be transformed to a state-transition model [32]. And this model can be further improved both in tightness and compactness by replace the original objective function and (23)(24) with the new objective function (25) and new constraint (26). Denoted this resulted model as 3P-Ti-ST [32].

## III. METHODOLOGY FOR MODELING IN THREE PERIODS

### A. Unit power variables

According to (8)-(9) and $u_{i,t} \in \{0,1\}$, we have $P_{i,t} \in 0 \cup [\underline{P}_i, \overline{P}_i]$. Then $P_{i,t}$ can be called as semi-continuous variable. Let [37], [12]

$$\tilde{P}_{i,t} = \frac{P_{i,t} - u_{i,t}\underline{P}_i}{(\overline{P}_i - \underline{P}_i)}, \quad (36)$$

then, we have $\tilde{P}_{i,t} \in [0,1]$ which represents the proportion, in $[\underline{P}_i, \overline{P}_i]$, of the power produced above $\underline{P}_i$ by unit $i$ at period $t$. Similarly, we let $\tilde{P}_{\text{up},i} = \frac{P_{\text{up},i}}{\overline{P}_i - \underline{P}_i}$, $\tilde{P}_{\text{down},i} = \frac{P_{\text{down},i}}{\overline{P}_i - \underline{P}_i}$, $\tilde{P}_{\text{start},i} = \frac{P_{\text{start},i} - \underline{P}_i}{\overline{P}_i - \underline{P}_i}$, $\tilde{P}_{\text{shut},i} = \frac{P_{\text{shut},i} - \underline{P}_i}{\overline{P}_i - \underline{P}_i}$. Then, $\tilde{P}_{\text{up},i} \in (0,1]$ represents the proportion of $P_{\text{up},i}$ in $[\underline{P}_i, \overline{P}_i]$. $\tilde{P}_{\text{down},i} \in (0,1]$, $\tilde{P}_{\text{start},i} \in [0,1]$, $\tilde{P}_{\text{shut},i} \in [0,1]$ have the similar means.

In actuality, [14] suggested use variable $P'_{i,t}$, recently has been put forward by [15], to represent the power produced above $\underline{P}_i$, which is more concise than our former transformation (36). However, with (36), the production cost $f_i(P_{i,t})$ can be transformed to a "tighter" form [12]

$$\tilde{f}_i(\tilde{P}_{i,t}) = \tilde{\alpha}_i u_{i,t} + \tilde{\beta}_i \tilde{P}_{i,t} + \tilde{\gamma}_i (\tilde{P}_{i,t})^2, \quad (37)$$

where $\tilde{\alpha}_i = \alpha_i + \beta_i \underline{P}_i + \gamma_i (\underline{P}_i)^2$, $\tilde{\beta}_i = (\overline{P}_i - \underline{P}_i)(\beta_i + 2\gamma_i \underline{P}_i)$, and $\tilde{\gamma}_i = \gamma_i(\overline{P}_i - \underline{P}_i)^2$.

Then the power balance constraint can be reformulated as:

$$\sum_{i=1}^{N}\left[\tilde{P}_{i,t}(\overline{P}_i - \underline{P}_i) + u_{i,t}\underline{P}_i\right] - P_{D,t} = 0. \quad (38)$$

### B. New state variables and logical constraints

In order to present our high dimension formulation for UC problem, we introduce more state variables for unit $i$ in periods $\{t-1, t, t+1\}$:

$o_{i,t}$: 1 if unit $i$ remains operational at time $t$ (i.e., $u_{i,t-1} = u_{i,t} = 1$), 0 otherwise.

$f_{i,t}$: 1 if unit $i$ remains off at time $t$, 0 otherwise.

$\mathcal{T}^1_{i,t} \sim \mathcal{T}^8_{i,t}$: There are eight combination for $u_{i,t}$ with periods

TABLE I
ILLUSTRATION FOR STATE VARIABLES AND UPPER BOUND OF GENERATION LIMITS

| $u_{t-1}$ | $u_t$ | $u_{t+1}$ | $o_t$ | $f_t$ | $s_t$ | $d_t$ | $o_{t+1}$ | $f_{t+1}$ | $s_{t+1}$ | $d_{t+1}$ | $\mathcal{T}^1_t$ | $\mathcal{T}^2_t$ | $\mathcal{T}^3_t$ | $\mathcal{T}^4_t$ | $\mathcal{T}^5_t$ | $\mathcal{T}^6_t$ | $\mathcal{T}^7_t$ | $\mathcal{T}^8_t$ | $\tilde{P}_t$ | UB($\tilde{P}_t$) | UB($\tilde{P}_{t-1}$) | UB($\tilde{P}_{t+1}$) |
|---|---|---|---|---|---|---|---|---|---|---|---|---|---|---|---|---|---|---|---|---|---|---|
| 0 | 0 | 0 | 0 | 1 | 0 | 0 | 0 | 1 | 0 | 0 | 1 | 0 | 0 | 0 | 0 | 0 | 0 | 0 | 0 | 0 | 0 | 0 |
| 0 | 0 | 1 | 0 | 1 | 0 | 0 | 0 | 0 | 1 | 0 | 0 | 1 | 0 | 0 | 0 | 0 | 0 | 0 | 0 | 0 | 0 | $\tilde{P}_{\text{start}}$ |
| 0 | 1 | 0 | 0 | 0 | 1 | 0 | 0 | 0 | 0 | 1 | 0 | 0 | 1 | 0 | 0 | 0 | 0 | 0 | $\tilde{P}_{i,t}$ | min($\tilde{P}_{\text{start}}, \tilde{P}_{\text{shut}}$) | 0 | 0 |
| 0 | 1 | 1 | 0 | 0 | 1 | 0 | 1 | 0 | 0 | 0 | 0 | 0 | 0 | 1 | 0 | 0 | 0 | 0 | $\tilde{P}_{i,t}$ | $\tilde{P}_{\text{start}}$ | 0 | $\tilde{P}_{\text{start}} + \tilde{P}_{\text{up}}; 1$ |
| 1 | 0 | 0 | 0 | 0 | 0 | 1 | 0 | 1 | 0 | 0 | 0 | 0 | 0 | 0 | 1 | 0 | 0 | 0 | 0 | 0 | $\tilde{P}_{\text{shut}}$ | 0 |
| 1 | 0 | 1 | 0 | 0 | 0 | 1 | 0 | 0 | 1 | 0 | 0 | 0 | 0 | 0 | 0 | 1 | 0 | 0 | 0 | 0 | $\tilde{P}_{\text{shut}}$ | $\tilde{P}_{\text{start}}$ |
| 1 | 1 | 0 | 1 | 0 | 0 | 0 | 0 | 0 | 0 | 1 | 0 | 0 | 0 | 0 | 0 | 0 | 1 | 0 | $\tilde{P}_{i,t}$ | $\tilde{P}_{\text{shut}}$ | $\tilde{P}_{\text{shut}} + \tilde{P}_{\text{down}}; 1$ | 0 |
| 1 | 1 | 1 | 1 | 0 | 0 | 0 | 1 | 0 | 0 | 0 | 0 | 0 | 0 | 0 | 0 | 0 | 0 | 1 | $\tilde{P}_{i,t}$ | 1 | 1 | 1 |

$\{t-1, t, t+1\}$. We use eights indicator binary variables $\mathcal{T}_{i,t}^1 \sim \mathcal{T}_{i,t}^8$ to indicate each status. For example, $\mathcal{T}_{i,t}^1 = 1$ if $u_{i,t-1} = u_{i,t} = u_{i,t+1} = 0$, $\mathcal{T}_{i,t}^1 = 0$ otherwise. All these state variables are illustrated in Table I.

We should note that, $o_{i,t}$ has been introduced in [32], and most of these state variables will not really appear in our final formulations. However, we introduce these state variables for two reasons: 1) facilitating our modeling; 2) illustrating the relationships of these state variables. Actually, it is not difficult to verify that the column vectors under $u_{t-1}, u_t, u_{t+1}, s_t$(or $d_t$), $s_{t+1}$(or $d_{t+1}$), $\mathcal{T}_t^2, \mathcal{T}_t^3$, and $e = [1; ...; 1]$ are linearly independent. These vectors can be viewed as base vectors and all the other state variables listed in Table I can be expressed as linear combination of these base vectors. For instance,

$$f_t = e - s_t - u_{t-1}, \quad (39)$$
$$\mathcal{T}_t^5 = s_t - u_t + u_{t-1} - s_{t+1} + \mathcal{T}_t^2 = d_t - s_{t+1} + \mathcal{T}_t^2, \quad (40)$$
$$\mathcal{T}_t^4 = s_t - \mathcal{T}_t^3, \quad (41)$$
$$\mathcal{T}_t^7 = s_{t+1} - u_{t+1} + u_t - \mathcal{T}_t^3 = d_{t+1} - \mathcal{T}_t^3, \quad (42)$$
$$\mathcal{T}_t^8 = -s_t - s_{t+1} + u_{t+1} + \mathcal{T}_t^3 = u_t - s_t - d_{t+1} + \mathcal{T}_t^3. \quad (43)$$

Actually, the logical constraints (4) and (5) also can be directly deduced based on these base vectors. We note that, either $s_t$ or $d_t$ can be equivalently chosen to be base vector, and we do it based on the compactness of the resulted formulations.

*C. The high dimension three-periods locally ideal model*

We use UB(·) to denote upper bound on "·", and RHS to represent right hand side. Consider the following inequality for generation limit,

$$P_{i,t} \leq \text{RHS expression}. \quad (44)$$

According to the physics meaning of generation limit, it is obvious that,

$$\text{UB}(P_{i,t}) = \begin{cases} \overline{P}_i & \text{if } u_{i,t} = 1 \\ 0 & \text{if } u_{i,t} = 0. \end{cases}$$

According to the UB$(P_{i,t})$, we let "RHS expression" of (44) be defined as $u_{i,t}\overline{P}_i$. Then we have $P_{i,t} \leq u_{i,t}\overline{P}_i$, i.e. (9) in this paper, the simplest generation limit.

According to the front analysis, one should try to compress the upper bound on left hand side (LHS) of (44). Other physics constraints can be considered in this procedure, and construct RHS expression with equal value with UB(LHS) for all possible statuses. Then the stronger possible inequality would be obtained.

When ramping limits in three periods are taking consideration, we list UB$(\tilde{P}_t)$ in Table I.

According to Table I, for $t \in \{2, ..., T-1\}$, we have

$$\text{UB}(P_{i,t}) = \begin{cases} 0 & \text{if } \mathcal{T}_t^1 = 1 \\ 0 & \text{if } \mathcal{T}_t^2 = 1 \\ \min(\tilde{P}_{\text{start}}, \tilde{P}_{\text{shut}}) & \text{if } \mathcal{T}_t^3 = 1 \\ \tilde{P}_{\text{start}} & \text{if } \mathcal{T}_t^4 = 1 \\ 0 & \text{if } \mathcal{T}_t^5 = 1 \\ 0 & \text{if } \mathcal{T}_t^6 = 1 \\ \tilde{P}_{\text{shut}} & \text{if } \mathcal{T}_t^7 = 1 \\ 1 & \text{if } \mathcal{T}_t^8 = 1 \end{cases}$$

i.e., $\text{UB}(\tilde{P}_t) = \mathcal{T}_t^3\{\min(\tilde{P}_{\text{start}}, \tilde{P}_{\text{shut}})\} + \mathcal{T}_t^4 \tilde{P}_{\text{start}} + \mathcal{T}_t^7 \tilde{P}_{\text{shut}} + \mathcal{T}_t^8$, then we obtain our upper bound limit for power of unit $i$, i.e.,

$$\tilde{P}_t \leq \mathcal{T}_t^3\{\min(\tilde{P}_{\text{start}}, \tilde{P}_{\text{shut}})\} + \mathcal{T}_t^4 \tilde{P}_{\text{start}} + \mathcal{T}_t^7 \tilde{P}_{\text{shut}} + \mathcal{T}_t^8 \quad (45)$$

According to linear relations (41)-(43) based on base vectors, we can transform (45) as

$$\tilde{P}_t \leq u_t - s_t(1 - \tilde{P}_{\text{start}}) - d_{t+1}(1 - \tilde{P}_{\text{shut}}) + \mathcal{T}_t^3\{1 - \max(\tilde{P}_{\text{start}}, \tilde{P}_{\text{shut}})\}, \quad (46)$$

Then, we obtain a strong valid inequality for unit power limit upper bound with additional state variables $\mathcal{T}_t^3$, and according to Table I, $\mathcal{T}_t^3$ can be determined by the following inequalities

$$\mathcal{T}_t^3 \geq s_t + d_{t+1} - u_t, \quad (47)$$
$$\mathcal{T}_t^3 \leq s_t, \quad (48)$$
$$\mathcal{T}_t^3 \leq d_{t+1} \quad (49)$$

We also list UB$(\tilde{P}_{t-1})$ and UB$(\tilde{P}_{t+1})$ in Table I. In these columns, we put the value of "UB(·)" on the right of a semicolon when only considering start-up and shut-down ramping limits, and put the value of "UB(·)" on the left of this semicolon when considering the whole ramping limits in three periods. And we only fill one value if there was no difference between these two cases.

In this paper, for the sake of convenience, we assume that $\tilde{P}_{\text{shut}} + \tilde{P}_{\text{down}} < 1$ and $\tilde{P}_{\text{start}} + \tilde{P}_{\text{up}} < 1$, otherwise, the operation ramping limits actually are absence. Actually, we only need to replace $\tilde{P}_{\text{shut}} + \tilde{P}_{\text{down}}$ and $\tilde{P}_{\text{start}} + \tilde{P}_{\text{up}}$ with $\min\{\tilde{P}_{\text{shut}} + \tilde{P}_{\text{down}}, 1\}$ and $\min\{\tilde{P}_{\text{start}} + \tilde{P}_{\text{up}}, 1\}$ respectively even if we are not sure the aforementioned assumptions. Similarly, we also assume $2\tilde{P}_{\text{up}} < 1$ and $2\tilde{P}_{\text{down}} < 1$.

Then, like the construction of (46), when the whole ramping limits in three periods are taking consideration, we obtain the following tight upper bounds,

$$\tilde{P}_{t-1} \leq \mathcal{T}_t^3(1 - \tilde{P}_{\text{down}} - \tilde{P}_{\text{shut}}) + u_{t-1} + d_t(\tilde{P}_{\text{shut}} - 1) + d_{t+1}(\tilde{P}_{\text{down}} + \tilde{P}_{\text{shut}} - 1) \quad (50)$$

$$\tilde{P}_{t+1} \leq \mathcal{T}_t^3(1 - \tilde{P}_{\text{start}} - \tilde{P}_{\text{up}}) + u_{t+1} + s_t(\tilde{P}_{\text{start}} + \tilde{P}_{\text{up}} - 1) + s_{t+1}(\tilde{P}_{\text{start}} - 1). \quad (51)$$

When only start-up and shut-down ramping limits are taking consideration, the following constraints, which are equivalent to (13) and (12) respectively, are obtained.

$$\tilde{P}_{t-1} \leq u_{t-1} - d_t(1 - \tilde{P}_{\text{shut}}) \quad (52)$$
$$\tilde{P}_{t+1} \leq u_{t+1} - s_{t+1}(1 - \tilde{P}_{\text{start}}). \quad (53)$$

For the sake of clarity, we neglected the history of the generators by removing the corresponding turn on/off constraints defined for $u_{i,t} = u_{i,0}, t \in [1, ..., U_i + L_i]$, but the following statements are correct for the case with history of the generators.

We use $\dim(\mathcal{P})$ to denote the dimension of a polytope $\mathcal{P}$. A valid inequality is a facet for $\mathcal{P}$ if and only if there are $\dim(\mathcal{P})$ affinely independent solutions in $\mathcal{P}$ that satisfy the inequality with equality [29]. Now we use this property to show the tightness of our high dimension generation limits.

Let $\mathcal{B}_t^R = \{u_t \in \mathcal{R}_+^3 \times [0,1]^8 | (4)(6) - (7)(46) - (51)\}$, where $u_t = (\tilde{P}_{t-1}, \tilde{P}_t, \tilde{P}_{t+1}, u_{t-1}, u_t, u_{t+1}, s_t, s_{t+1}, d_t, d_{t+1}, \mathcal{T}_t^3)$. $\mathcal{B}_t^I = \mathcal{B}_t^R \cap (\mathcal{R}_+^3 \times \{0,1\}^8)$.

TABLE II
AFFINELY INDEPENDENT POINTS FOR CONSTRAINTS OF $\mathcal{B}_t^I$ and $\mathcal{C}_t^I$

| Point | $\tilde{P}_{t-1}$ | $\tilde{P}_t$ | $\tilde{P}_{t+1}$ | $u_{t-1}$ | $u_t$ | $u_{t+1}$ | $s_t$ | $s_{t+1}$ | $d_t$ | $d_{t+1}$ | $\mathcal{T}_t^3$ |
|---|---|---|---|---|---|---|---|---|---|---|---|
| $u_t^1$ | 0 | $\tilde{P}_{\text{start}}$ | $\tilde{P}_{\text{start}} + \tilde{P}_{\text{up}}$ | 0 | 1 | 1 | 1 | 0 | 0 | 0 | 0 |
| $u_t^2$ | 0 | 0 | $\tilde{P}_{\text{start}}$ | 0 | 0 | 1 | 0 | 1 | 0 | 0 | 0 |
| $u_t^3$ | $\tilde{P}_{\text{shut}}$ | 0 | 0 | 1 | 0 | 0 | 0 | 0 | 1 | 0 | 0 |
| $u_t^4$ | $\tilde{P}_{\text{start}} + \tilde{P}_{\text{down}}$ | $\tilde{P}_{\text{shut}}$ | 0 | 1 | 1 | 0 | 0 | 0 | 0 | 1 | 0 |
| $u_t^5$ | 0 | $\min(\tilde{P}_{\text{start}}, \tilde{P}_{\text{shut}})$ | 0 | 0 | 1 | 0 | 1 | 0 | 0 | 1 | 1 |
| $u_t^6$ | $\tilde{p}_1$ | $\tilde{p}_2$ | $\tilde{p}_3$ | 1 | 1 | 1 | 0 | 0 | 0 | 0 | 0 |
| $u_t^7$ | $\tilde{p}_4$ | $\tilde{p}_5$ | $\tilde{p}_6$ | 1 | 1 | 1 | 0 | 0 | 0 | 0 | 0 |
| $u_t^8$ | $\tilde{p}_7$ | $\tilde{p}_8$ | $\tilde{p}_9$ | 1 | 1 | 1 | 0 | 0 | 0 | 0 | 0 |
| $u_t^9$ | 0 | 0 | 0 | 0 | 0 | 0 | 0 | 0 | 0 | 0 | 0 |
| $u_t^{10}$ | 0 | 0 | 0 | 1 | 1 | 1 | 0 | 0 | 0 | 0 | 0 |
| $u_t^{11}$ | 0 | 0 | $\tilde{P}_{\text{start}}$ | 1 | 0 | 1 | 0 | 1 | 1 | 0 | 0 |
| $u_t^{12}$ | 0 | $\tilde{P}_{\text{start}}$ | 0 | 0 | 1 | 1 | 1 | 0 | 0 | 0 | 0 |
| $u_t^{13}$ | 0 | $\tilde{P}_{\text{shut}} - \varepsilon$ | 0 | 1 | 1 | 0 | 0 | 0 | 0 | 1 | 0 |
| $u_t^{14}$ | 0 | $\tilde{P}_{\text{start}}$ | $[\tilde{P}_{\text{start}} - \tilde{P}_{\text{down}}]^+$ | 0 | 1 | 1 | 1 | 0 | 0 | 0 | 0 |
| $u_t^{15}$ | $[\tilde{P}_{\text{shut}} - \tilde{P}_{\text{up}}]^+$ | $\tilde{P}_{\text{shut}} - \varepsilon$ | 0 | 1 | 1 | 0 | 0 | 0 | 0 | 1 | 0 |

$\varepsilon > 0$ is a sufficiently small positive number.

**Theorem 1**: Given any $t \in \{2, \ldots, T-1\}$, when $i \in \mathcal{J}^1$, $\dim(\text{conv}(\mathcal{B}_t^I)) = 9$ and the inequalities (46)-(49)(50)(51) describe facets of the polytope $\text{conv}(\mathcal{B}_t^I)$ respectively, and for $i \in \mathcal{J}^{\geq 2}$, $\dim(\text{conv}(\mathcal{B}_t^I)) = 8$ and the inequalities (46)(50)(51) with $\mathcal{T}_t^3 = 0$ still describe facets of the polytope $\text{conv}(\mathcal{B}_t^I)$ respectively. And $\text{conv}(\mathcal{B}_t^I) = \mathcal{B}_t^R$.

Proof. At first, we proof the conclusion for $i \in \mathcal{J}^1$. Let $\mathcal{A}_{\text{UB}}$ denotes $\begin{bmatrix} \tilde{p}_1 & \tilde{p}_2 & \tilde{p}_3 \\ \tilde{p}_4 & \tilde{p}_5 & \tilde{p}_6 \\ \tilde{p}_7 & \tilde{p}_8 & \tilde{p}_9 \end{bmatrix}$ in Table II. Observing the linear equation (4) in $\mathcal{B}_t^I$, we know that $\dim(\text{conv}(\mathcal{B}_t^I))$ can at most be 9. This is indeed true because that the 10 points $(u_t^k, k = 1, \ldots, 10)$ of $\mathcal{B}_t^I$ listed in Table II are affinely-independent when $\mathcal{A}_{\text{UB}} = \begin{bmatrix} 1 & 1 & 0 \\ 0 & 1 & 0 \\ 1 & 1 & 1 \end{bmatrix}$. In this case, the 9 points $(u_t^k, k = 1, \ldots, 9)$ are tight for inequality (46). Then (46) describes a facet of $\text{conv}(\mathcal{B}_t^I)$.

The 9 points $(u_t^k, k = 2, \ldots, 10)$ are tight for inequality (48). Then (48) describes a facet of $\text{conv}(\mathcal{B}_t^I)$. The 9 points $(u_t^k, k = 1,2,3,5, \ldots, 10)$ are tight for inequality (49). Then (49) describes a facet of $\text{conv}(\mathcal{B}_t^I)$. The 9 affinely-independent points $(u_t^k, k = 1, \ldots, 5, 9, 11, \ldots, 13)$ are tight for inequality (47). Then (47) describes a facet of $\text{conv}(\mathcal{B}_t^I)$.

Similarly, (50) describes a facet of $\text{conv}(\mathcal{B}_t^I)$ because that the affinely-independent 9 points $(u_t^k, k = 1, \ldots, 9)$ are tight for inequality (50) when $\mathcal{A}_{\text{UB}} = \begin{bmatrix} 1 & 0 & 0 \\ 1 & 1 & 0 \\ 1 & 0 & 1 \end{bmatrix}$. Finally, (51) describes a facet of $\text{conv}(\mathcal{B}_t^I)$ because that the affinely-independent 9 points $(u_t^k, k = 1, \ldots, 9)$ are tight for inequality (51) when $\mathcal{A}_{\text{UB}} = \begin{bmatrix} 1 & 0 & 1 \\ 0 & 1 & 1 \\ 0 & 0 & 1 \end{bmatrix}$.

Now, we will prove that $\text{conv}(\mathcal{B}_t^I) = \mathcal{B}_t^R$.

Projecting out $d$ from (47)(48)(49) by using (4), we have
$$\mathcal{T}_t^3 \geq s_t + s_{t+1} - u_{t+1}, \quad (54)$$
$$\mathcal{T}_t^3 \leq s_t, \quad (55)$$
$$\mathcal{T}_t^3 \leq s_{t+1} - u_{t+1} - u_t \quad (56)$$

We use $\mathcal{B}_t^T$ to denote the set $\{y_t \in [0,1]^6 | (54) - (56), (6) - (7)\}$, where $y_t = (u_{t-1}, u_t, u_{t+1}, s_t, s_{t+1}, \mathcal{T}_t^3)$, and (7) is rewritten by projecting out $d$ accordingly. It is not difficult to verify that $\mathcal{B}_t^T$ is an integral polyhedron because that the coefficient matrix of $\mathcal{B}_{\mathcal{T},t}^3$ is totally unimodular (by using pivoting operation) [29]. We give the full details as follows.

$$\mathcal{B}_t^T = \left\{ y_t \in [0,1]^6 \middle| \begin{array}{c} \mathcal{T}_t^3 \geq s_t + s_{t+1} - u_{t+1} \\ \mathcal{T}_t^3 \leq s_t \\ s_t \leq u_t \\ s_t \leq 1 - u_{t-1} \\ s_{t+1} \leq u_{t+1} \\ s_{t+1} \leq 1 - u_t \end{array} \right\}$$

The coefficient matrix of $\mathcal{B}_t^T$ is:
$$B = \begin{bmatrix} 0 & 0 & -1 & 1 & 1 & -1 \\ 0 & 0 & 0 & -1 & 0 & 1 \\ 0 & -1 & 0 & 1 & 0 & 0 \\ 1 & 0 & 0 & 1 & 0 & 0 \\ 0 & 0 & -1 & 0 & 1 & 0 \\ 0 & 1 & 0 & 0 & 1 & 0 \end{bmatrix}$$

After the axis operation of row and column for $B$, we get:
$$B' = \begin{bmatrix} 0 & 0 & 0 & 1 & 0 & -1 \\ 0 & 0 & 0 & 0 & 0 & 0 \\ 0 & 0 & 0 & 0 & 1 & 1 \\ 1 & 0 & 0 & 0 & 0 & 0 \\ 0 & 0 & 1 & 0 & 0 & 0 \\ 0 & 1 & 0 & -1 & 0 & 0 \end{bmatrix}$$

According to corollary 3.4 in [38], $B'$ is totally unimodular. Let $a = c = [0; 0; 0; 0; 0; 0]$, $b = [1; 1; 1; 1; 1; 1]$, $d = [0; 0; 0; 1; 0; 1]$, we know that the vertices of $\mathcal{B}_t^T$ all are integer according to corollary 3.2 in [38].

Finally, according to the Lemma 4 in [17], it is not difficult to prove $\text{conv}(\mathcal{B}_t^I) = \mathcal{B}_t^R$.

For $i \in \mathcal{J}^{\geq 2}$, it is obvious that $u_t^5$ is an infeasible status and $\mathcal{T}_t^3 = 0$, so (47)-(49) are redundant which can be removed from $\mathcal{B}_t^I$. Then a similar proof can be conducted. ∎

*Remark 1*: When we replace (50)(51) in definition of $\mathcal{B}_t^R$

TABLE III
ILLUSTRATION FOR UPPER BOUND OF RAMPING LIMITS

| $u_{t-1}$ | $u_t$ | $u_{t+1}$ | UB$(\tilde{P}_t - \tilde{P}_{t-1})$ | UB$(\tilde{P}_t - \tilde{P}_{t+1})$ | UB$(\tilde{P}_{t+1} - \tilde{P}_t)$ | UB$(\tilde{P}_{t-1} - \tilde{P}_t)$ | UB$(\tilde{P}_{t+1} - \tilde{P}_{t-1})$ | UB$(\tilde{P}_{t-1} - \tilde{P}_{t+1})$ |
|---|---|---|---|---|---|---|---|---|
| 0 | 0 | 0 | 0 | 0 | 0 | 0 | 0 | 0 |
| 0 | 0 | 1 | 0 | 0 | $\tilde{P}_{\text{start}}$ | 0 | $\tilde{P}_{\text{start}}$ | 0 |
| 0 | 1 | 0 | $\min(\tilde{P}_{\text{start}}, \tilde{P}_{\text{shut}})$ | $\min(\tilde{P}_{\text{start}}, \tilde{P}_{\text{shut}})$ | 0 | 0 | 0 | 0 |
| 0 | 1 | 1 | $\tilde{P}_{\text{start}}$ | $\min(\tilde{P}_{\text{start}}, \tilde{P}_{\text{down}})$ | $\tilde{P}_{\text{up}}$ | 0 | $\tilde{P}_{\text{start}} + \tilde{P}_{\text{up}}$ | 0 |
| 1 | 0 | 0 | 0 | 0 | 0 | $\tilde{P}_{\text{shut}}$ | 0 | $\tilde{P}_{\text{shut}}$ |
| 1 | 0 | 1 | 0 | 0 | $\tilde{P}_{\text{start}}$ | $\tilde{P}_{\text{shut}}$ | $\tilde{P}_{\text{start}}$ | $\tilde{P}_{\text{shut}}$ |
| 1 | 1 | 0 | $\min(\tilde{P}_{\text{up}}, \tilde{P}_{\text{shut}})$ | $\tilde{P}_{\text{shut}}$ | 0 | $\tilde{P}_{\text{down}}$ | 0 | $\tilde{P}_{\text{shut}} + \tilde{P}_{\text{down}}$ |
| 1 | 1 | 1 | $\tilde{P}_{\text{up}}$ | $\tilde{P}_{\text{down}}$ | $\tilde{P}_{\text{up}}$ | $\tilde{P}_{\text{down}}$ | $2\tilde{P}_{\text{up}}$ | $2\tilde{P}_{\text{down}}$ |

with (52)(53) and use $\widetilde{\mathcal{B}}_t^R$ to denote the resulted set, a conclusion similar to Theorem 1 can be obtain for $\widetilde{\mathcal{B}}_t^R$. With some minor revisions for $u_t^1$ and $u_t^4$, the proof of this new conclusion is nearly as same as the proof of Theorem 1. $\widetilde{\mathcal{B}}_T^R \supseteq \mathcal{B}_T^R$.

Let $\mathcal{B}_T^R = \{(P, u, s, d, \mathcal{T}^3) \in \mathcal{R}_+^T \times [0,1]^{4T-4} | (4)(6)(7)(46) - (51)\}$ and $\mathcal{B}_T^I = \mathcal{B}_T^R \cap (\mathcal{R}_+^T \times \{0,1\}^{4T-4})$. Following a similar but more complicated proof, the Theorem 1 can be generalized as:

**Corollary 1**: when $i \in \mathcal{I}^1$, $\dim(\text{conv}(\mathcal{B}_T^I)) = 4T - 3$ and the inequalities (46)-(51) describe facets of the polytope $\text{conv}(\mathcal{B}_T^I)$ respectively, and for $i \in \mathcal{I}^{\geq 2}$, $\dim(\text{conv}(\mathcal{B}_T^I)) = 3T - 1$ and the inequalities (46)(50)(51) with $\mathcal{T}_t^3 = 0$ still describe facets of the polytope $\text{conv}(\mathcal{B}_T^I)$ respectively. And $\text{conv}(\mathcal{B}_T^I) = \mathcal{B}_T^R$.

*Remark 2*: Similar to remark 1, when we replace (50)(51) in definition of $\mathcal{B}_T^R$ with (52)(53) and use $\widetilde{\mathcal{B}}_T^R$ to denote the resulted set, a conclusion similar to Corollary 1 can be obtain for $\widetilde{\mathcal{B}}_T^R$. And we note that, for $\widetilde{\mathcal{B}}_T^R$ with $i \in \mathcal{I}^{\geq 2}$, an equivalent conclusion of corollary 1 was given in [17].

For power limits inequalities, we have the following conclusions:

1) (48)(36)(46)⇒(14)⇒(12); (46)≳(14)≳(12);
2) (49)(36)(46)⇒(15)⇒(13); (46)≳(15)≳(13);
3) For $i \in \mathcal{I}^{\geq 2}$, $\mathcal{T}_t^3$ vanishes identically, then (36)(46)⇒(11) ⇒(10)⇒(9), (46)≂(11)≳(10)≳(9).
4) (50)≳(52), (51)≳(53), (46)≳(52)(53). And (50) is *strictly tighter* than (52) when $\tilde{P}_{\text{shut}} + \tilde{P}_{\text{down}} < 1$; (51) is *strictly tighter* than (53) when $\tilde{P}_{\text{start}} + \tilde{P}_{\text{up}} < 1$.

These conclusions show that our proposed power generation upper bound inequality imply the several state-of-art inequality listed in Section II.B.3. Here we note that, [18][19][20], and [35] provided power generation upper bound inequalities with multiple periods. However, our formulations are tighter than these inequalities when only three periods are considered.

Similar to analysis of upper bound limit for power of unit $i$, when startup and shutdown ramping limits are taking consideration, we list UB$(\tilde{P}_t - \tilde{P}_{t-1})$ for ramp-up limit in Table III.

According to Table III, for $t \in \{2, ..., T-1\}$, we obtain our ramp-up inequality for unit $i$.

$$\tilde{P}_t - \tilde{P}_{t-1} \leq \mathcal{T}_t^3 \min(\tilde{P}_{\text{start}}, \tilde{P}_{\text{shut}}) + \mathcal{T}_t^4 \tilde{P}_{\text{start}} + \mathcal{T}_t^7 \min(\tilde{P}_{\text{up}}, \tilde{P}_{\text{shut}}) + \mathcal{T}_t^8 \tilde{P}_{\text{up}} \quad (57)$$

With (41)(42)(43)(5), then (57) can be further transformed as

$$\tilde{P}_t - \tilde{P}_{t-1} \leq \mathcal{T}_t^3 \left([\tilde{P}_{\text{up}} - \tilde{P}_{\text{shut}}]^+ - [\tilde{P}_{\text{start}} - \tilde{P}_{\text{shut}}]^+\right) +$$
$$u_t \tilde{P}_{\text{up}} + s_t(\tilde{P}_{\text{start}} - \tilde{P}_{\text{up}}) - d_{t+1}[\tilde{P}_{\text{up}} - \tilde{P}_{\text{shut}}]^+ \quad (58)$$

When bounding ramping over three periods, according to the penultimate column in Table III, we have

$$\tilde{P}_{t+1} - \tilde{P}_{t-1} \leq \mathcal{T}_t^3(\tilde{P}_{\text{up}} - \tilde{P}_{\text{start}}) + u_{t+1} 2\tilde{P}_{\text{up}} + s_t(\tilde{P}_{\text{start}} - \tilde{P}_{\text{up}}) + s_{t+1}(\tilde{P}_{\text{start}} - 2\tilde{P}_{\text{up}}) \quad (59)$$

In addition, according to Table III, we have

$$\tilde{P}_{t+1} - \tilde{P}_t \leq u_{t+1}\tilde{P}_{\text{up}} + s_{t+1}(\tilde{P}_{\text{start}} - \tilde{P}_{\text{up}}) \quad (60)$$

Similarly, according to Table III, we have our ramp-down constraint as follows,

$$\tilde{P}_t - \tilde{P}_{t+1} \leq \mathcal{T}_t^3\left([\tilde{P}_{\text{down}} - \tilde{P}_{\text{start}}]^+ - [\tilde{P}_{\text{shut}} - \tilde{P}_{\text{start}}]^+\right) + u_t \tilde{P}_{\text{down}} - s_t[\tilde{P}_{\text{down}} - \tilde{P}_{\text{start}}]^+ + d_{t+1}(\tilde{P}_{\text{shut}} - \tilde{P}_{\text{down}}) \quad (61)$$

$$\tilde{P}_{t-1} - \tilde{P}_{t+1} \leq \mathcal{T}_t^3(\tilde{P}_{\text{down}} - \tilde{P}_{\text{shut}}) + u_{t-1} 2\tilde{P}_{\text{down}} + d_t(\tilde{P}_{\text{shut}} - 2\tilde{P}_{\text{down}}) + d_{t+1}(\tilde{P}_{\text{shut}} - \tilde{P}_{\text{down}}) \quad (62)$$

$$\tilde{P}_{t-1} - \tilde{P}_t \leq u_{t-1}\tilde{P}_{\text{down}} + d_t(\tilde{P}_{\text{shut}} - \tilde{P}_{\text{down}}) \quad (63)$$

Let $\mathcal{C}_t^R = \mathcal{B}_t^R \cap \{u_t|(58) - (63)\}$, $\mathcal{C}_t^I = \mathcal{C}_t^R \cap (\mathcal{R}_+^3 \times \{0,1\}^8)$. It is not very hard to get the following Theorem:

**Theorem 2**: Given any $t \in \{2, ..., T - 1\}$, when $i \in \mathcal{I}^1$, $\dim(\text{conv}(\mathcal{C}_t^I)) = 9$ and the inequalities (46)-(51), (58)-(63) describe facets of the polytope $\text{conv}(\mathcal{C}_t^I)$ respectively, and for $i \in \mathcal{I}^{\geq 2}$, $\dim(\text{conv}(\mathcal{C}_t^I)) = 8$ and the inequalities (46), (58)-(63) describe facets of the polytope $\text{conv}(\mathcal{B}_t^I)$ respectively.

Proof. At first, we proof the conclusion for $i \in \mathcal{I}^1$. Observing the linear equation (4) in $\mathcal{C}_t^I$, we know that $\dim(\text{conv}(\mathcal{C}_t^I))$ can at most be 9. This is indeed true because that the 10 points $(u_t^k, k = 1, ..., 10)$ of $\mathcal{C}_I^t$ listed in Table II are affinely-independent when $\mathcal{A}_{\text{UB}} = \begin{bmatrix} 1 & 1 & 1 - \tilde{P}_{\text{down}} \\ 1 - \tilde{P}_{\text{up}} & 1 & 1 - \tilde{P}_{\text{down}} \\ 1 & 1 & 1 \end{bmatrix}$.

In this case, the 9 points $(u_t^k, k = 1, ..., 9$ in Table II) are tight for inequality (46). Then (46) describes a facet of $\text{conv}(\mathcal{C}_t^I)$. The 9 points $(u_t^k, k = 2, ..., 10$ in Table II) are tight for inequality (48). Then (48) describes a facet of $\text{conv}(\mathcal{C}_t^I)$. The 9 points $(u_t^k, k = 1,2,3,5, ..., 10$ in Table II) are tight for inequality (49). Then (49) describes a facet of $\text{conv}(\mathcal{C}_t^I)$. The 9 affinely-independent points $(u_t^k, k = 1, ..., 5, 9, 11, 14, 15$ in Table II) are tight for inequality (47). Then (47) describes a facet of $\text{conv}(\mathcal{C}_t^I)$.

Similarly, (50) describes a facet of $\text{conv}(\mathcal{C}_t^I)$ because that the affinely-independent 9 points $(u_t^k, k = 1, ..., 9$ in Table II) are tight for inequality (50) when $\mathcal{A}_{\text{UB}} = \begin{bmatrix} 1 & 1 - \tilde{P}_{\text{down}} & 1 - \tilde{P}_{\text{down}} \\ 1 & 1 & 1 - \tilde{P}_{\text{down}} \\ 1 & 1 & 1 \end{bmatrix}$.

Finally, (51) describes a facet of $\text{conv}(\mathcal{C}_t^I)$ because that the affinely-independent 9 points ($u_t^k, k = 1, \dots, 9$) are tight for inequality (51) when $\mathcal{A}_{\text{UB}} = \begin{bmatrix} 1 & 1 & 1 \\ 1 - \tilde{P}_{\text{up}} & 1 & 1 \\ 1 - \tilde{P}_{\text{up}} & 1 - \tilde{P}_{\text{up}} & 1 \end{bmatrix}$.

At last, (58)-(63) describe facets of $\text{conv}(\mathcal{C}_t^I)$ respectively because that the affinely-independent points listed in Table A.I of Appendix A are tight for these inequalities respectively.

For $i \in \mathcal{J}^{\geq 2}$, it is obvious that $u_t^5$ is an infeasible status and $\mathcal{T}_t^3 = 0$, so (47)-(49) are redundant which can be removed from $\mathcal{C}_t^I$. Then a similar proof can be conducted. ∎

We note that here, as far as we know, more constraints (including up and down time, generation limits, and up and down ramping) are included in $\mathcal{C}_t^I$ simultaneously than other references. For instance, no reference considers up and down ramp constrains simultaneously, up and down time constraints have not been included in [18].

*Remark 3*: Different to Theorem 1, when we replace $\mathcal{B}_t^R$ with $\widetilde{\mathcal{B}}_t^R$ in Theorem 2, weak constraints (52) and (53) are not the facets of $\text{conv}(\mathcal{C}_t^I)$ because that $\mathcal{C}_t^I$ is invariable whatever we chose $\mathcal{B}_t^R$ or $\widetilde{\mathcal{B}}_t^R$.

In addition, we should point it out that the corollary of Theorem 2 similar to Corollary 1 is not true also.

For ramp limits inequalities, we have the following conclusions:

1) For $i \in \mathcal{J}^{\geq 2}$, $\mathcal{T}_t^3$ vanishes identically, then when $\tilde{P}_{\text{up}} > \tilde{P}_{\text{shut}}$, we have (36)(58)⇒(20), and (58) is *strictly tighter* than (20) when $\tilde{P}_{\text{up}} < \tilde{P}_{\text{shut}}$ (considering the state of $\mathcal{T}_t^7 = 1$);

2) Similarly, for $i \in \mathcal{J}^{\geq 2}$, when $\tilde{P}_{\text{down}} > \tilde{P}_{\text{start}}$, we have (36)(61)⇒(22), and (61) is *strictly tighter* than (22) when $\tilde{P}_{\text{down},i} < \tilde{P}_{\text{start},i}$;

3) For $i \in \mathcal{J}^{\geq 2} \cap \underline{\mathcal{J}}^{\geq 2}$, (36)(59)⇒(21);

4) (36)(58)⇒(18)(16), and (58) is strictly tighter than (18) when $\tilde{P}_{\text{shut}} < \max\{\tilde{P}_{\text{start}}, \tilde{P}_{\text{up}}\}$. Actually, with the $\tilde{P}_t$ variables, constraint (18) and (16) can be easily constructed according to Table II and only considering periods $t - 1$ and $t$;

5) (36)(61)⇒(19)(17), and (61) is strictly tighter than (19) when $\tilde{P}_{\text{start}} < \max\{\tilde{P}_{\text{shut}}, \tilde{P}_{\text{down}}\}$. Similarly, with the $\tilde{P}_t$ variables, constraint (19) and (17) can be easily constructed according to Table III.

Now, we present our tight and compact MIQP UC formulation in *high dimension* space (3P-HD):

$$\min F_C = \sum_{i=1}^N \sum_{t=1}^T [\tilde{f}_i(\tilde{P}_{i,t}) + C_{\text{hot},i} s_{i,t} + \tilde{S}_{i,t}]$$

$$\text{s.t.} \begin{cases} (26)(4)(32)(6)(7) \\ (47)(48)(49) \\ (46)(50)(51) \\ (58)(59)(61)(62) \\ (60) \text{ for } t = T - 1, (63) \text{ for } t = 2 \\ (38)(3) \\ \mathcal{T}_{i,t}^3 = 0 \text{ for } i \in \mathcal{J}^{\geq 2} \\ (u_{i,t}, s_{i,t}, d_{i,t}, \mathcal{T}_{i,t}^3) \in \{0,1\}^4, (\tilde{P}_{i,t}, \tilde{S}_{i,t}) \in \mathcal{R}_+^2. \end{cases} \quad (64)$$

In the following content of this subsection, we will project the 3P-HD model onto the set of variables for classic 3-bin UC.

We note that here, for $i \in \mathcal{J}^{\geq 2}$, $\mathcal{T}_t^3$ vanishes identically, then (58) is equivalent to

$$\tilde{P}_t - \tilde{P}_{t-1} \leq u_t \tilde{P}_{\text{up}} + s_t (\tilde{P}_{\text{start}} - \tilde{P}_{\text{up}}) - d_{t+1} [\tilde{P}_{\text{up}} - \tilde{P}_{\text{shut}}]^+ \quad (65)$$

Let $\varrho_i^1 = [\tilde{P}_{i,\text{up}} - \tilde{P}_{i,\text{shut}}]^+ - [\tilde{P}_{i,\text{start}} - \tilde{P}_{i,\text{shut}}]^+$, then, for $i \in \mathcal{J}^1$ and $\varrho_i^1 \geq 0$, considering (48)(49), (58) is equivalent to (because $\mathcal{T}_t^3 = \min\{s_t, d_{t+1}\}$) the following two inequalities:

$$\tilde{P}_t - \tilde{P}_{t-1} \leq s_t (\min\{\tilde{P}_{\text{start}}, \tilde{P}_{\text{shut}}\} - \min\{\tilde{P}_{\text{up}}, \tilde{P}_{\text{shut}}\}) + u_t \tilde{P}_{\text{up}} - d_{t+1} [\tilde{P}_{\text{up}} - \tilde{P}_{\text{shut}}]^+, \quad (66)$$

$$\tilde{P}_t - \tilde{P}_{t-1} \leq s_t (\tilde{P}_{\text{start}} - \tilde{P}_{\text{up}}) + u_t \tilde{P}_{\text{up}} - d_{t+1} [\tilde{P}_{\text{start}} - \tilde{P}_{\text{shut}}]^+. \quad (67)$$

Here we note that (58)≳(66)(67).

For $i \in \mathcal{J}^1$ and $\varrho_i^1 < 0$, considering (47), (58) is equivalent to two inequalities, i.e. (65) and the following inequality (because $\mathcal{T}_t^3 = \max\{s_t + d_{t+1} - u_t, 0\}$)

$$\tilde{P}_t - \tilde{P}_{t-1} \leq -d_{t+1} [\tilde{P}_{\text{start}} - \tilde{P}_{\text{shut}}]^+ + u_t (\min\{\tilde{P}_{\text{up}}, \tilde{P}_{\text{shut}}\} + [\tilde{P}_{\text{start}} - \tilde{P}_{\text{shut}}]^+) + s_t (\min\{\tilde{P}_{\text{start}}, \tilde{P}_{\text{shut}}\} - \min\{\tilde{P}_{\text{up}}, \tilde{P}_{\text{shut}}\}) \quad (68)$$

That is, (58) can be transformed to be (65) or (66)(67) or (68)(65) after eliminating the variable $\mathcal{T}_t^3$.

Similarly, after eliminating $\mathcal{T}_t^3$ from (50), (51), (59), (61), and (62), we can obtain the corresponding constraints (73)-(85) (listed in Appendix B). We should note that, most of the resulted inequalities after eliminating $\mathcal{T}_t^3$ from (50), (51), (58), (59), (61), and (62), for all we know, have not been found on known publications.

In addition, $\mathcal{T}_t^3$ also can be eliminated from (46) *without losing tightness* (Because (11)(14)(15) have been proved to be facets and can form a convex hull together with the minimum up-times and down-times constraint [17]) as shown in appendix (86)(87), and the equivalent inequalities have also been obtained in [15] and [17] ((11)(14)(15) in this paper).

Then our 3P-HD model can be projected into traditional UC space, and we denote this model as 3P-HD-Pr.

$$\min F_C = \sum_{i=1}^N \sum_{t=1}^T [\tilde{f}_i(\tilde{P}_{i,t}) + C_{\text{hot},i} s_{i,t} + \tilde{S}_{i,t}]$$

$$\text{s.t.} \begin{cases} (26)(4)(32)(6)(7) \\ (86) \text{ for } i \in \mathcal{J}^{\geq 2}; (87) \text{ for } i \in \mathcal{J}^1, \\ (73)(75) \text{ for } i \in \mathcal{J}^{\geq 2}; (74)^2 (76)^2 \text{ for } i \in \mathcal{J}^1; \\ (52) \text{ for } t = 2, (53) \text{ for } t = T - 1, \\ (65)(77)(80)(83) \text{ for } i \in \mathcal{J}^{\geq 2} \\ (66)(67) \text{ for } i \in \mathcal{J}^1 \text{ and } \varrho_i^1 \geq 0 \\ (68)(65) \text{ for } i \in \mathcal{J}^1 \text{ and } \varrho_i^1 < 0 \\ (78)(79)(81)(82)(84)(85) \\ (60) \text{ for } t = T - 1, (63) \text{ for } t = 2 \\ (38)(3) \\ (u_{i,t}, s_{i,t}, d_{i,t}) \in \{0,1\}^3, (\tilde{P}_{i,t}, \tilde{S}_{i,t}) \in \mathcal{R}_+^2, \end{cases} \quad (69)$$

At last, according to the descriptions of our models and the other four models reviewed in the previous section, we have the following relationship for these models in tightness,

2P-Co ≾ 2P-Ti ≾ 3P-Ti ≾ 3P-Ti-ST ≾ 3P-HD-Pr ≾ 3P-HD.

$$\tag{70}$$

## D. MILP Approximations

Solving an MIQP formulation for UC as an MILP using linear approximation is very popular because the significant improvements in off-the-shelf MILP solvers.

Assume that $L$ is a given parameter, let $p_{i,l} = \underline{P}_i + l\,(\overline{P}_i - \underline{P}_i)/L$ and $l = 0,1,2,\ldots,L$. For 2P-Co, 2P-Ti, and 3P-Ti, after replacing $\gamma_i(P_{i,t})^2$ in the objective function with a corresponding new variable $z_{i,t}$ and adding the following linear constraints to the formulation [11]

$$z_{i,t} \geq 2\gamma_i p_{i,l} P_{i,t} - \gamma_i(p_{i,l})^2, \tag{71}$$

we obtain the MILP UC models which approximate the original MIQP models. For 3P-T-ST, (5) should be use to replace $u_{i,t}$ with $o_{i,t} + s_{i,t}$.

Similarly, let $\tilde{p}_{i,l} = l/L$. Replace $\tilde{\gamma}_i(\tilde{P}_{i,t})^2$ for 3P-HD, and 3P-HD-Pr models with $z_{i,t}$ and add the following constraints

$$z_{i,t} \geq 2\tilde{\gamma}_i \tilde{p}_{i,l} \tilde{P}_{i,t} - \tilde{\gamma}_i(\tilde{p}_{i,l})^2, \tag{72}$$

and then, the corresponding MILP approximations have been formed.

## IV. NUMERICAL RESULTS AND ANALYSIS

In this section, we present some numerical results to test the efficiency and effectiveness [34] of the proposed 3P-HD and 3P-HD-Pr MIP formulations. 51 randomly generated realistic instances with units running from 10 to 1080 for a time span of 24 hours are used in our experiments. There are 5 different data values each for 10-, 20-, 50-, 75-, 100- and 150-unit systems, 12 data values for 200-unit systems [39], one data value each for 560-, 700-, 880-, 900-, 980-, 1000-, 1040- and 1080-unit system, two data values for 1020-unit systems [19], respectively. For the rest of this paper, all the results reported for 10 to 200-unit systems are the averaged results of systems with the same number of units, unless otherwise noted. Table IV shows the No. and the number of units for each test instance. The machine on which we perform all of our computations is Dell XPS8930 with 16 GB of RAM and Intel i7-8700 3.2GHz CPU, running MS-Windows 10 and Matlab2016b. We use CPLEX 12.7.1 to solve MILP and MIQP. Note that the time limit for the solver is 3600 s. The source code of our experiments can be found and download at [https://github.com/linfengYang/High-dimension-UC-formulation].

TABLE IV NUMBER OF UNITS IN EACH PROBLEM INSTANCE

| No. | The number of unit | No. | The number of unit |
|---|---|---|---|
| 1 | 10 | 10 | 880 |
| 2 | 20 | 11 | 900 |
| 3 | 50 | 12 | 980 |
| 4 | 75 | 13 | 1000 |
| 5 | 100 | 14 | 1020 |
| 6 | 150 | 15 | 1020 |
| 7 | 200 | 16 | 1040 |
| 8 | 560 | 17 | 1080 |
| 9 | 700 | - | - |

## A. Comparison of 6 MIQP formulations

We are now in a position to show the tightness of our proposed 3P-HD and 3P-HD-Pr MIQP formulation. The tightness of an MIP formulation is defined as the deference of optimal values between continuous relaxation and the original MIP problem, so the tightness often can be measured by using relative integrality gap [29]. The relative integrality gap of a given MIP can be defined as $(Z_{\text{MIP}} - Z_{\text{CR}})/Z_{\text{MIP}}$. In this expression, $Z_{\text{CR}}$ represents the optimal value of continuous relaxation of original MIP problem and $Z_{\text{MIP}}$ represents the optimal value of this MIP. However, in practice, the MIP problems are hard to solve for real optimal but only to a preset optimal tolerance. To make the comparison to be fair, we use the same $Z_{\text{MIP}}$ for six MIQPs, that is the best integer solution which was found by using CPLEX with accuracy 0.5% among aforementioned six MIQPs. And we denote the initial relative integrality gap calculated in this way as "iGap". Then, the difference in iGap will accurately depict the difference in tightness of MIP's initial continuous relaxations.

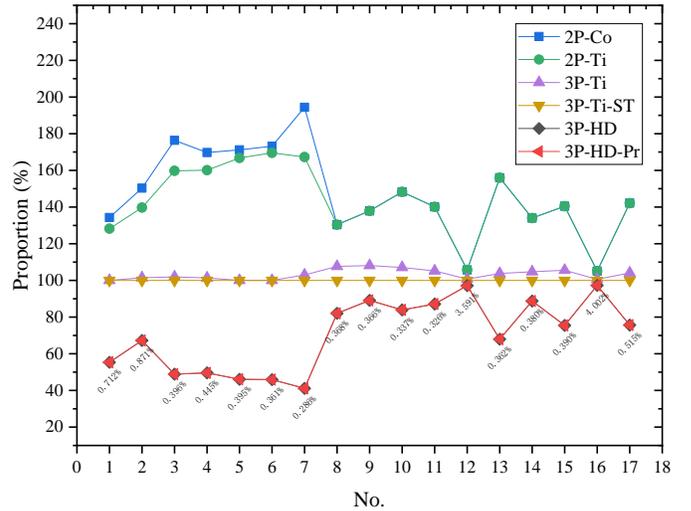

Fig.1. Comparsion of the tightness of the six MIQP formulations in terms of iGap.

Fig. 1 gives iGap for all the formulations in comparison with 3P-Ti-ST's (using ratios) for all test systems, where the ratio for 3P-Ti-ST model itself always represents 100%. And the detailed values of iGap for 3P-HD model have been labeled in the figure as well. In this figure, "No." represents the No. of each problem instance. As can be seen in Fig. 1, our 3P-HD and 3P-HD-Pr models have the same tightness and are tighter than the other models. Meanwhile, 3P-Ti-ST is tighter than 3P-Ti for most cases. But the improvement effect of 3P-Ti-ST in tightness is not significant for the test instances of 10 to 200-unit systems. The main reason is that 3P-Ti-ST mainly improves the tightness by changing the startup cost constraint (as analyzed in Section II), and the startup cost of these test instances only accounts for a small part of the total cost. However, the improvement effect of our two models is very significant for these test instances. The reason of these results is that the ramping re-strictions of these instances are randomly generated more realistic and complicated than other 8-unit based test instances.

After directly solving the six MIQP models using CPLEX to 0.5% optimality within 3600 s, we report the CPU time in Fig. 2 and detailed results in Table V. It should be noted that when fast response times are required in a given operating environment, 0.5% is a very popular tolerance to specify for CPLEX to ensure that UC MIP models are solved within a reasonable execution time with reasonably satisfying suboptimal solutions.

In order to better express the difference of results, in Fig. 2, we introduce the relative time, which can be defined as $(C_{\text{time\_MIP}} - C_{\text{time\_REF}})/C_{\text{time\_REF}}$. In this expression, $C_{\text{time\_MIP}}$ represents the CPU time of original MIP problem and $C_{\text{time\_REF}}$ represents the CPU time of the reference model. We denote the relative time calculated in this way as "rTime". And here we still choose 3P-Ti-ST model as the reference model. The detailed values of CPU time for 3P-Ti-ST model have been labeled in the figure as well.

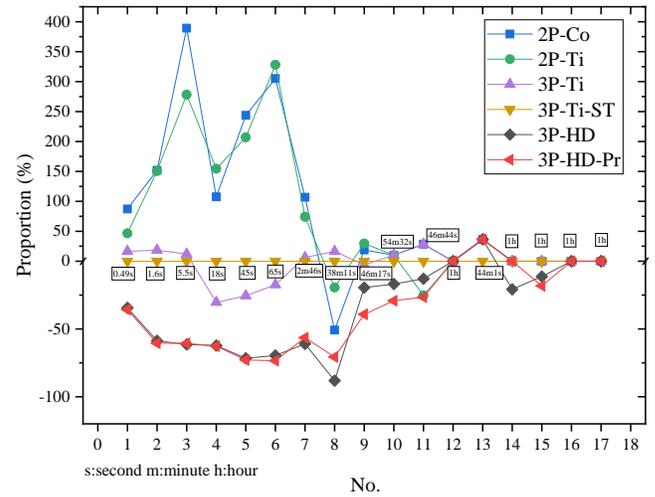

Fig.2. Comparsion of the six MIQP formulations in terms of rTime.

Fig.2 gives rTime for all the formulations in comparison with

TABLE V COMPARISON OF SIX ORIGINAL MIQPS IN THE COMPUTATIONAL PERFORMANCE

| No. | 2P-Co | | | | | | | 2P-Ti | | | | | |
|---|---|---|---|---|---|---|---|---|---|---|---|---|---|
| | $C_{\text{time}}$(s) | Gap(%) | Nodes | Cuts | $N_{\text{it\_all}}$ | $I_u$(%) | $I_{\text{all}}$(%) | $C_{\text{time}}$(s) | Gap(%) | Nodes | Cuts | $N_{\text{it\_all}}$ | $I_u$(%) | $I_{\text{all}}$(%) |
| 1 | 0.92 | 0.49 | 9 | 280 | 3086 | 48.67 | 69.72 | 0.72 | 0.39 | 11 | 199 | 4204 | 49.83 | 72.22 |
| 2 | 3.97 | 0.48 | 213 | 696 | 131517 | 50.75 | 70.76 | 3.95 | 0.47 | 347 | 559 | 110468 | 52.5 | 73.5 |
| 3 | 26.69 | 0.42 | 320 | 1687 | 22614 | 52.88 | 71.63 | 20.63 | 0.46 | 303 | 1239 | 46697 | 55.08 | 74.65 |
| 4 | 38.14 | 0.45 | 112 | 2450 | 66222 | 53.66 | 72.23 | 46.77 | 0.45 | 228 | 2053 | 69486 | 55.08 | 74.34 |
| 5 | 153.22 | 0.47 | 638 | 3202 | 125519 | 52.78 | 72.39 | 136.8 | 0.47 | 809 | 2824 | 132934 | 53.95 | 73.58 |
| 6 | 265.78 | 0.39 | 606 | 4692 | 152243 | 55.28 | 73.93 | 280.88 | 0.41 | 1045 | 4163 | 122424 | 56.37 | 74.99 |
| 7 | 345.03 | 0.38 | 726 | 6508 | 161195 | 55.58 | 72.94 | 290.39 | 0.41 | 692 | 4430 | 141626 | 57.41 | 76.76 |
| 8 | 1127.34 | 0.48 | 1911 | 3582 | 467901 | 59.43 | 80.44 | 1847.20 | 0.50 | 2186 | 1604 | 533848 | 59.49 | 80.45 |
| 9 | 3300.58 | 0.35 | 4821 | 3993 | 512584 | 48.67 | 74.05 | 3600.84 | 2.68 | 3013 | 2917 | 431239 | 48.67 | 74.08 |
| 10 | 3601.77 | 3.11 | 4681 | 7475 | 429301 | 42.74 | 70.10 | 3601.52 | 4.16 | 4182 | 4595 | 448433 | 42.78 | 70.10 |
| 11 | 3601.82 | 4.45 | 4411 | 5500 | 251552 | 56.30 | 77.68 | 2096.59 | 0.25 | **3581** | 3715 | 206693 | 56.30 | 77.68 |
| 12 | 3601.69 | 3.4 | 3911 | 9324 | 358670 | 42.98 | 69.88 | 3600.31 | 4.09 | **3103** | 10753 | 350274 | 42.98 | 69.90 |
| 13 | 3601.45 | 5.51 | 2100 | 12327 | 333068 | 62 | 77.92 | 3601.59 | 6.19 | 3464 | 10322 | 332148 | 62 | 77.92 |
| 14 | 3601.41 | 0.28 | 4737 | 7653 | 299511 | 56.70 | 79.07 | 3601.94 | 4.45 | **2800** | 3616 | 260032 | 56.73 | 79.06 |
| 15 | 3601.55 | 4.77 | 4022 | 9407 | 520674 | 55.24 | 75.58 | 3601.05 | 5.80 | **2237** | 6267 | 458228 | 55.24 | 75.58 |
| 16 | 3601.23 | **3.81** | 4214 | 10627 | 285637 | 44.95 | 70.29 | 3601.86 | 4.49 | 4342 | 6679 | 472281 | 44.95 | 70.29 |
| 17 | 3600.53 | 5.84 | **2763** | 10832 | 367300 | 61.74 | 76.92 | 3601.48 | 6.16 | 3995 | 8047 | 515373 | 61.77 | 76.85 |

| No. | 3P-Ti | | | | | | | 3P-Ti-ST | | | | | | |
|---|---|---|---|---|---|---|---|---|---|---|---|---|---|---|
| | $C_{\text{time}}$(s) | Gap(%) | Nodes | Cuts | $N_{\text{it\_all}}$ | $I_u$(%) | $I_{\text{all}}$(%) | $C_{\text{time}}$(s) | Gap(%) | Nodes | Cuts | $N_{\text{it\_all}}$ | $I_u$(%) | $I_{\text{all}}$(%) |
| 1 | 0.57 | 0.41 | 1 | 124 | 1793 | 55.5 | 78.69 | 0.49 | 0.42 | 1 | 115 | 2596 | 54.25 | 78.22 |
| 2 | 2.16 | 0.45 | 5 | 323 | 6456 | 59.69 | 80.12 | 1.58 | 0.45 | 31 | 264 | 7237 | 58.75 | 79.89 |
| 3 | 6.11 | 0.42 | 5 | 853 | 4228 | 61.42 | 80.54 | 5.46 | 0.42 | 1 | 710 | 7869 | 60.47 | 80.36 |
| 4 | 12.79 | 0.43 | 2 | 1197 | 5781 | 63.69 | 81.60 | 18.38 | 0.45 | 79 | 1152 | 10974 | 63 | 81.46 |
| 5 | 33.22 | 0.42 | 95 | 1663 | 26797 | 64.08 | 82.05 | 44.58 | 0.40 | 221 | 1571 | 28130 | 63.02 | 81.71 |
| 6 | 54.24 | 0.42 | 173 | 2504 | 26817 | 65.51 | 82.70 | 65.61 | 0.34 | 212 | 2284 | 50190 | 64.47 | 82.37 |
| 7 | 177.34 | 0.30 | 542 | 3084 | 89717 | 63.10 | 81.55 | 166.81 | 0.36 | 633 | 2693 | 115847 | 62.42 | 81.48 |
| 8 | 2667.92 | 0.39 | 3544 | 886 | 276139 | 59.97 | 80.16 | 2291.23 | 0.42 | 2700 | **576** | 369869 | 59.97 | 80.01 |
| 9 | 2704.58 | 0.33 | 4685 | 1072 | 489614 | 47.5 | 73.52 | 2777.74 | 0.24 | 3897 | 779 | 687062 | 48.01 | 73.78 |
| 10 | 3601.92 | 4.00 | 3946 | 3852 | 255100 | 53.6 | 74.93 | 3272.17 | 0.15 | **3732** | 3029 | 241625 | 53.48 | 74.05 |
| 11 | 3601.56 | 5.02 | 4704 | 1767 | 215304 | 71.73 | 84.17 | 2804.64 | 0.12 | 3678 | 1925 | 246501 | 69.84 | 83.07 |
| 12 | 3601.86 | 4.06 | 3802 | 4821 | 407659 | 58.57 | 77.06 | 3602.25 | 4.012 | 4636 | 3593 | 250046 | 58.09 | 75.87 |
| 13 | 3601.42 | 0.15 | **1493** | 4908 | 328731 | 64.17 | 79.67 | **2641.70** | **0.09** | 3385 | 3999 | 237980 | 64 | 79.69 |
| 14 | 3602.23 | 4.31 | 4518 | 1788 | 194981 | 60.62 | 78.39 | 3601.88 | **0.11** | 4100 | 1222 | 207189 | 60.21 | 77.78 |
| 15 | 3602.55 | 5.63 | 3103 | 4181 | 347086 | 65.48 | 80.18 | 3601.19 | 5.65 | 5284 | 3991 | 235980 | 64.78 | 79.20 |
| 16 | 3602.08 | 4.36 | 4293 | 4203 | 502358 | 56.00 | 75.37 | 3601.11 | 4.36 | 5822 | 3857 | 310221 | 54.68 | 74.44 |
| 17 | 3601.81 | 5.90 | 3841 | 5430 | 277535 | 61.73 | 77.83 | 3601.77 | 5.93 | 3421 | 5464 | 202271 | 63.58 | 79.12 |

| No. | 3P-HD | | | | | | | 3P-HD-Pr | | | | | | |
|---|---|---|---|---|---|---|---|---|---|---|---|---|---|---|
| | $C_{\text{time}}$(s) | Gap(%) | Nodes | Cuts | $N_{\text{it\_all}}$ | $I_u$(%) | $I_{\text{all}}$(%) | $C_{\text{time}}$(s) | Gap(%) | Nodes | Cuts | $N_{\text{it\_all}}$ | $I_u$(%) | $I_{\text{all}}$(%) |
| 1 | **0.32** | **0.32** | **0** | 84 | 793 | **67.58** | **88.5** | **0.32** | 0.34 | **0** | 87 | **740** | 67.08 | 84.53 |
| 2 | **0.65** | **0.38** | **0** | 85 | 1063 | **75** | **90.90** | 0.62 | 0.40 | **0** | 110 | 1185 | 74.29 | 87.76 |
| 3 | **2.1** | **0.38** | **0** | 224 | 1518 | **79.7** | **92.31** | 2.15 | **0.33** | **0** | 206 | **1455** | 78.72 | 89.49 |
| 4 | **6.94** | **0.26** | **0** | 424 | **2721** | **80.58** | **92.49** | 6.86 | 0.33 | **0** | 421 | 2843 | 79.59 | 89.71 |
| 5 | 12.62 | 0.32 | **0** | **539** | 2651 | **81.03** | **92.81** | 12.11 | **0.21** | **0** | 550 | 3015 | 80.78 | 90.32 |
| 6 | 19.97 | **0.21** | 49 | 829 | 4029 | **83.13** | **93.60** | 17.42 | 0.27 | 47 | 797 | 3686 | 82.79 | 91.34 |
| 7 | **65.04** | **0.16** | 132 | 1228 | 27757 | **81.29** | **92.92** | 72.69 | **0.11** | 146 | 1304 | 21391 | 80.91 | 90.56 |
| 8 | **272.52** | **0.12** | **0** | 2042 | 149089 | 59.97 | **83.07** | 671.23 | 0.23 | 1503 | 2051 | **78004** | 60.30 | 79.86 |
| 9 | 2237.89 | **0.08** | 4076 | 310 | 313715 | 47.38 | **78.39** | 1687.83 | 0.22 | 3002 | 310 | 515665 | 49.39 | 73.89 |
| 10 | 2720.08 | **0.09** | 4317 | 656 | 145941 | 56.16 | 80.85 | 2316.94 | **0.05** | 4859 | 660 | 151536 | 57.86 | 76.15 |
| 11 | 2440.81 | **0.03** | 4187 | 955 | 169590 | 72.59 | **88.03** | 2062.5 | 0.11 | 4029 | **801** | 180310 | **72.87** | 85.10 |
| 12 | 3601.75 | 4.01 | 3620 | 2474 | **210903** | **62.24** | **84.57** | 3601.83 | 4.01 | 4080 | **2447** | 251375 | 61.81 | 79.42 |
| 13 | 3601.55 | 6.04 | 2017 | **4151** | 183777 | 73.08 | **88.04** | 3602.06 | 6.04 | 2530 | 4529 | 161013 | **73.50** | 84.66 |
| 14 | **2855.09** | **0.15** | 2985 | 15 | 179159 | 66.26 | **85.29** | 3601.45 | 0.14 | 3757 | **13** | 165289 | 67.54 | 82.32 |
| 15 | 3193.81 | **0.04** | 4270 | 2156 | **123510** | 68.46 | **86.19** | 2945.36 | 0.03 | 3906 | 2142 | 160852 | **70.02** | 83.00 |
| 16 | 3600.44 | 4.36 | 4988 | 1421 | **214299** | 58.81 | **83.29** | 3600.38 | 4.36 | 5039 | **612** | 226795 | 59.61 | 78.57 |
| 17 | 3601.97 | 5.92 | 4299 | **4189** | 137686 | 67.58 | **88.5** | 3602.11 | **0.15** | 4392 | 4236 | 161700 | **68.38** | 82.58 |

3P-Ti-ST's (using ratios) for all test systems. According to the definition of rTime, we know that the rTime for 3P-Ti-ST model is zero. It can be seen from the Fig.2, the results of our 3P-HD and 3P-HD-Pr models are better than those of the other four models for most instances. For 1020-, 1040- and 1080-unit system, the values of rTime for six models are roughly the same. This is because the problem size is too large for CPLEX to get a satisfactory solution in the limited time (3600 s).

In Table V, "$C_{\text{time}}$" represents the execution time, "Gap" represents the relative integrality gap reported by CPLEX at the end of solution process, "$N_{\text{it\_all}}$" represents the total number of iterations required to solve the node relaxations during the current optimization, "Cuts" represents the number of cuts generated by CPLEX to tighten the model, and "Nodes" represents the number of nodes visited while the solver is building the enumeration tree during the solution process, "$I_u$" represents the proportion of integers in $u_{i,t}$ for the solution of continuous relaxation, "$I_{\text{all}}$" represents the proportion of integers in the solution of continuous relaxation.

According to Table V, our 3P-HD and 3P-HD-Pr models clearly outperforms the other four formulations in terms of the computation time. These results are in full agreement with comparison of rTime presented in Fig. 2. Moreover, we observe that the number of nodes required to solve the 3P-HD and 3P-HD-Pr formulations is far less than that for other four formulations. The reason is that markedly fewer nodes need to be visited during the branch-and-cut procedure. Furthermore, Table V shows that the values of "$I_u$" and "$I_{\text{all}}$" for the 3P-HD and 3P-HD-Pr models are much larger than those of the other four models, especially 3P-HD. This finding further proves that the 3P-HD and 3P-HD-Pr models are tighter than the other four models. These experimental results show that the 3P-HD and 3P-HD-Pr models have better performance than the other four models.

### B. Comparison of 6 MILP formulations

According to Section III.D, the six MIQP formulations can be approximated as six corresponding MILP formulations, with $L = 4$ for our experiments. Similar to Fig.1, Fig.2 and Table V,

TABLE VI COMPARISON OF SIX MILPS IN THE COMPUTATIONAL PERFORMANCE

| No. | 2P-Co | | | | | | | 2P-Ti | | | | | | |
|---|---|---|---|---|---|---|---|---|---|---|---|---|---|---|
| | $C_{\text{time}}$(s) | Gap(%) | Nodes | Cuts | $N_{\text{it\_all}}$ | $I_u$(%) | $I_{\text{all}}$(%) | $C_{\text{time}}$(s) | Gap(%) | Nodes | Cuts | $N_{\text{it\_all}}$ | $I_u$(%) | $I_{\text{all}}$(%) |
| 1 | 0.53 | 0.35 | 0 | 399 | 2494 | 50.17 | 70.69 | 0.48 | 0.39 | 0 | 317 | 2282 | 50 | 72.75 |
| 2 | 1.16 | 0.37 | 0 | 863 | 5298 | 54.54 | 73.40 | 0.99 | 0.45 | 0 | 678 | 4707 | 55.13 | 75.29 |
| 3 | 4.62 | 0.42 | 0 | 2507 | 6844 | 53.75 | 73.09 | 3.02 | 0.42 | 0 | 1747 | 5038 | 55.55 | 75.63 |
| 4 | 8.81 | 0.42 | 0 | 3731 | 9398 | 55.82 | 74.14 | 8.82 | 0.39 | 0 | 2990 | 8138 | 56.53 | 75.88 |
| 5 | 16.67 | 0.37 | 0 | 5093 | 13428 | 55.21 | 74.47 | 14.61 | 0.43 | 0 | 4255 | 11318 | 55.81 | 75.38 |
| 6 | 21.37 | **0.26** | 0 | 6780 | 20613 | 55.87 | 75.06 | 19.25 | 0.28 | 0 | 6681 | 22691 | 56.39 | 75.87 |
| 7 | 41.12 | 0.22 | 13 | 10061 | 39433 | 56.01 | 74.11 | 31.96 | 0.30 | 13 | 7370 | 22410 | 57.15 | 77.20 |
| 8 | 324.14 | 0.31 | 0 | 7454 | 97676 | 47.17 | 72.22 | 425.73 | 0.41 | 20 | 7016 | 132041 | 29.32 | 66.27 |
| 9 | 613.83 | 0.28 | 16 | 9973 | 245630 | 14.52 | 59.09 | 1125.44 | 0.30 | 291 | **8262** | 1138471 | 22.26 | 62.18 |
| 10 | 2835.20 | 0.16 | 232 | 18777 | 904035 | 21.40 | 61.02 | 2782.41 | 0.16 | 224 | 15058 | 1619764 | 41.29 | 67.65 |
| 11 | 322.09 | 0.22 | 0 | 9932 | 76153 | 60.83 | 78.55 | 1450.53 | 0.20 | 382 | 8537 | 916257 | 42.78 | 72.53 |
| 12 | 3291.52 | 0.44 | 161 | 20890 | 847170 | 37.84 | 66.50 | 3600.81 | 49.07 | 146 | 17257 | 903114 | 32.74 | 64.80 |
| 13 | 3601.05 | 38.53 | 852 | 18671 | 1015905 | 63.67 | 77.56 | 3600.81 | 49.41 | 1935 | 18455 | 1720746 | 63.33 | 77.44 |
| 14 | 1202.25 | 0.14 | 173 | 13405 | 979203 | 46.41 | 72.93 | 2659.48 | 0.23 | 388 | 8798 | 967837 | 54.74 | 75.71 |
| 15 | 3601.13 | 38.95 | 192 | 19184 | 1409278 | 60.46 | 76.80 | 3601.06 | 49.47 | 294 | 16439 | 1070891 | 60.46 | 76.80 |
| 16 | 3600.94 | 39.02 | 119 | 24242 | 2091648 | 26.28 | 61.70 | 3600.98 | 38.75 | 272 | 18374 | 2997365 | 40.71 | 66.51 |
| 17 | 3600.75 | 0.23 | 287 | 21146 | 1445081 | 60.03 | 75.72 | 3102.52 | 0.29 | 384 | 17629 | 2063445 | 53.86 | 73.61 |

| No. | 3P-Ti | | | | | | | 3P-Ti-ST | | | | | | |
|---|---|---|---|---|---|---|---|---|---|---|---|---|---|---|
| | $C_{\text{time}}$(s) | Gap(%) | Nodes | Cuts | $N_{\text{it\_all}}$ | $I_u$(%) | $I_{\text{all}}$(%) | $C_{\text{time}}$(s) | Gap(%) | Nodes | Cuts | $N_{\text{it\_all}}$ | $I_u$(%) | $I_{\text{all}}$(%) |
| 1 | 0.28 | 0.33 | 0 | 233 | 2118 | 58 | 79.75 | 0.37 | 0.39 | 0 | 374 | 2456 | 56.25 | 79.11 |
| 2 | 0.76 | 0.43 | 0 | 455 | 3220 | 59.83 | 80.47 | 0.87 | 0.38 | 0 | 614 | **2087** | 58.708 | 80.10 |
| 3 | 2.19 | 0.45 | 0 | 1232 | 2984 | 62.93 | 81.48 | 2.55 | 0.46 | 0 | 1433 | 3012 | 61.683 | 81.2 |
| 4 | 4.01 | 0.46 | 0 | 2029 | 3866 | 64.73 | 82.47 | 4.22 | 0.40 | 0 | 2250 | 4246 | 63.967 | 82.29 |
| 5 | 5.24 | 0.38 | 0 | 2902 | 6770 | 65.24 | 83.09 | 5.77 | 0.39 | 0 | 3550 | 6527 | 64.092 | 82.71 |
| 6 | 9.54 | 0.37 | 0 | 4636 | 8246 | 65.78 | 83.35 | 9.74 | 0.39 | 0 | 5316 | 8321 | 64.578 | 82.94 |
| 7 | 14.70 | 0.33 | 0 | 4843 | 12187 | 64.29 | 82.35 | 20.38 | **0.19** | 0 | 5884 | 13220 | 63.052 | 82.11 |
| 8 | 347.22 | 0.41 | 0 | **5846** | 63803 | 22.47 | 64.34 | 364.83 | 0.08 | 0 | 7473 | 74797 | 59.8214 | 77.33 |
| 9 | 894.70 | 0.30 | 129 | 8554 | **45345** | 22.86 | 62.82 | 495.92 | 0.17 | 0 | 11023 | 74150 | 27.5 | 64.05 |
| 10 | 466.75 | 0.26 | 0 | **10982** | 43581 | 44.13 | 71.24 | 406.09 | 0.16 | 0 | 12626 | 53220 | 46.4015 | 71.65 |
| 11 | 246.69 | 0.20 | 0 | 5418 | 35244 | 68.89 | 81.98 | 1279.92 | **0.03** | 125 | 10400 | 52909 | 68.3333 | 81.70 |
| 12 | 2302.14 | 0.24 | 149 | 13876 | 445592 | 50.60 | 73.75 | 526.05 | 0.12 | 0 | **11975** | 63222 | 54.3367 | 74.29 |
| 13 | 1644.81 | 0.25 | 492 | 14735 | 875480 | 62.5 | 79.11 | 1203.38 | 0.11 | 657 | 13191 | 606717 | 57.75 | 77.61 |
| 14 | **101.52** | 0.20 | 0 | **3555** | 538183 | 61.11 | 77.51 | 274.73 | 0.25 | 0 | 12870 | **72785** | 60.2124 | 77.10 |
| 15 | 587.09 | 0.23 | 0 | **10936** | 55618 | 60.21 | 78.10 | 338.16 | 0.12 | 0 | **10936** | 59977 | 38.7255 | 69.72 |
| 16 | 791.47 | 0.26 | 0 | 14311 | 70836 | 44.07 | 70.70 | 506.48 | 0.1 | 0 | **12110** | 68838 | 52.6442 | 73.16 |
| 17 | 2978.86 | 0.1 | 332 | **17190** | 634577 | 55.48 | 75.95 | 3601.08 | 0.15 | 397 | 17399 | 822980 | 61.5741 | 78.40 |

| No. | 3P-HD | | | | | | | 3P-HD-Pr | | | | | | |
|---|---|---|---|---|---|---|---|---|---|---|---|---|---|---|
| | $C_{\text{time}}$(s) | Gap(%) | Nodes | Cuts | $N_{\text{it\_all}}$ | $I_u$(%) | $I_{\text{all}}$(%) | $C_{\text{time}}$(s) | Gap(%) | Nodes | Cuts | $N_{\text{it\_all}}$ | $I_u$(%) | $I_{\text{all}}$(%) |
| 1 | **0.23** | 0.29 | 0 | **172** | 1697 | **67.58** | **88.5** | 0.27 | **0.25** | 0 | 272 | **1671** | 67.25 | 84.53 |
| 2 | **0.42** | **0.27** | 0 | **303** | 2874 | 75 | **90.90** | 0.48 | 0.39 | 0 | 322 | 2956 | **75.51** | 88.28 |
| 3 | **1.26** | 0.38 | 0 | **592** | 1068 | **79.7** | **92.30** | 1.31 | 0.39 | 0 | 629 | 1092 | 79.57 | 89.79 |
| 4 | 2.12 | **0.35** | 0 | 778 | 1553 | **80.58** | **92.49** | **1.83** | 0.39 | 0 | 634 | **1287** | 80.43 | 89.99 |
| 5 | 2.29 | 0.37 | 0 | 1007 | 2052 | 81.03 | **92.81** | **2.16** | 0.41 | 0 | 853 | **1576** | **81.05** | 90.42 |
| 6 | **3.87** | 0.36 | 0 | **1132** | 2432 | **83.13** | **93.60** | 3.92 | 0.40 | 0 | 1389 | 2909 | 83.09 | 91.45 |
| 7 | 6.03 | 0.36 | 0 | **1250** | 2627 | 81.29 | **92.92** | 6.34 | 0.38 | 0 | 1539 | 3109 | 81.13 | 90.63 |
| 8 | 195.70 | **0.04** | 0 | 11004 | 88516 | 59.97 | 83.07 | 120.08 | 0.40 | 0 | 8849 | 90998 | 59.97 | 77.43 |
| 9 | 338.84 | **0.17** | 0 | 14893 | 78237 | **47.38** | **78.39** | 104.52 | 0.22 | 0 | 10602 | 63378 | **47.38** | 71.19 |
| 10 | **228.72** | **0.09** | 0 | 14971 | 73597 | **56.16** | **80.85** | 833.67 | 0.12 | 122 | 17543 | 366250 | 55.78 | 74.49 |
| 11 | 175.59 | 0.06 | 0 | 10985 | 64047 | **72.59** | **88.03** | 158.92 | 0.06 | 0 | 11386 | 60752 | **72.59** | 84.14 |
| 12 | **321.70** | **0.05** | 0 | 19393 | **48994** | **62.24** | **84.57** | 339 | **0.05** | 0 | 19223 | 52339 | **62.24** | 79.34 |
| 13 | 1000.69 | 0.07 | 190 | 11264 | 670070 | **73.08** | **88.04** | 226.25 | **0.05** | 0 | **9912** | 35861 | 72.67 | 83.92 |
| 14 | 261.91 | 0.15 | 0 | 15793 | 97686 | **66.26** | **85.29** | 322.30 | **0.12** | 0 | 15710 | 449029 | **66.26** | 80.50 |
| 15 | 1624.90 | **0.04** | 233 | 15318 | 661056 | 68.46 | **86.19** | 1291.13 | **0.04** | 213 | 15916 | 634878 | 69.28 | 82.05 |
| 16 | 304.03 | **0.04** | 0 | 22032 | **43166** | **58.81** | **83.29** | **269.47** | 0.05 | 0 | 22426 | 45147 | **58.81** | 77.86 |
| 17 | 471.06 | 0.10 | 0 | 20692 | 88134 | **68.36** | **86.88** | 386.31 | **0.01** | 0 | 21325 | **74374** | **68.36** | 82.51 |

a comparison of the tightness for the six MILP formulations is shown in Fig.3, a comparison of rTime for the six MILP formulations is presented in Fig.4, and Table VI further shows the comparison of the six MILP formulations in details.

As can be seen in Fig. 3, our 3P-HD and 3P-HD-Pr models have the same tightness and they are tighter than the other models. And the detailed values of iGap for 3P-HD model have been labeled in Fig. 3. 3P-Ti-ST is tighter than 3P-Ti for most cases. It can be seen in Fig.4, the results of our 3P-HD and 3P-HD-Pr models are better than those of the other four models for most instances. In a word, Fig.3 and Fig.4 demonstrate that our 3P-HD and 3P-HD-Pr models are tighter and more effective than the other four models. The detailed values of CPU time for 3P-Ti-ST model have been labeled in Fig.4.

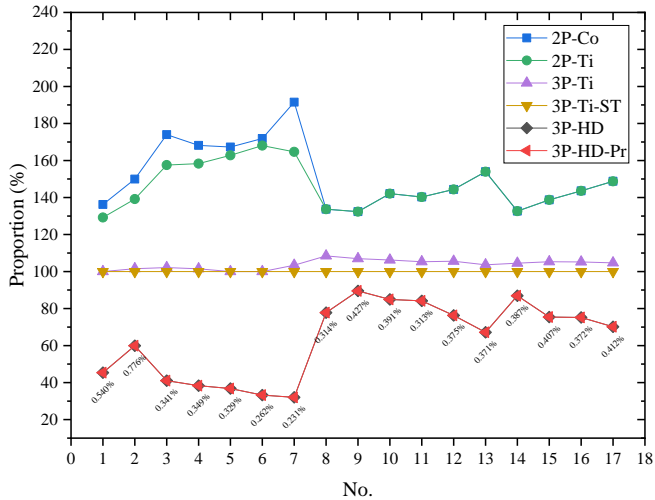

Fig.3. Comparsion of the tightness of the six MILP formulations in terms of iGap.

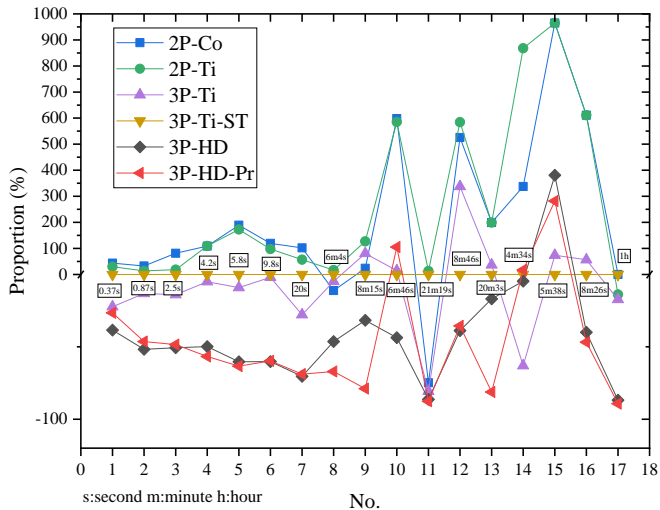

Fig.4. Comparsion of the six MILP formulations in terms of rTime.

As shown in Table VI, similar to MIQPs, the proposed MILP models are significantly better than the other four MILP models for most instances. For our models, the CPU times are often less, the total numbers of nodes required to solve are often fewer, the values of "$I_u$" and "$I_{all}$" are always higher, these test results indicate that the performance of the proposed MILP models is significantly better than the other four models.

## V. Conclusion

In this study, by introducing more auxiliary state variables, we developed a tight and compact high dimension UC model and presented the detailed modeling procedure. This model is locally ideal, i.e., we give the novel generations limits and production ramping limits which are the cornerstone constraints (facets) forming the convex hull description of a single unit in 3 periods. This model in high dimension space can be projected onto variables set of the classic 3-bin UC. And which deduces a new 3-bin formulation nearly as tight as the proposed high dimension model. We compare the performance of our new formulations with four benchmark formulations. Our proposed models perform the best due to their outstanding tightness without losing compactness too much.

APPENDIX A.

TABLE A.I

(1) AFFINELY INDEPENDENT POINTS FOR CONSTRAINTS (58)

| Point | $\tilde{P}_{t-1}$ | $\tilde{P}_t$ | $\tilde{P}_{t+1}$ | $u_{t-1}$ | $u_t$ | $u_{t+1}$ | $s_t$ | $s_{t+1}$ | $d_t$ | $d_{t+1}$ | $\mathcal{T}_t^3$ |
|---|---|---|---|---|---|---|---|---|---|---|---|
| $u_t^1$ | 0 | $\tilde{P}_{\text{start}}$ | $\tilde{P}_{\text{start}} + \tilde{P}_{\text{up}}$ | 0 | 1 | 1 | 1 | 0 | 0 | 0 | 0 |
| $u_t^2$ | 0 | 0 | $\tilde{P}_{\text{start}}$ | 0 | 0 | 1 | 0 | 1 | 0 | 0 | 0 |
| $u_t^3$ | 0 | 0 | 0 | 1 | 0 | 0 | 0 | 0 | 1 | 0 | 0 |
| $u_t^4$ | 0 | $\min(\tilde{P}_{\text{up}}, \tilde{P}_{\text{shut}})$ | 0 | 1 | 1 | 0 | 0 | 0 | 0 | 1 | 0 |
| $u_t^5$ | 0 | $\min(\tilde{P}_{\text{start}}, \tilde{P}_{\text{shut}})$ | 0 | 0 | 1 | 0 | 1 | 0 | 0 | 1 | 1 |
| $u_t^6$ | 0 | $\tilde{P}_{\text{up}}$ | $-\varepsilon + \tilde{P}_{\text{up}}$ | 1 | 1 | 1 | 0 | 0 | 0 | 0 | 0 |
| $u_t^7$ | $\varepsilon$ | $\varepsilon + \tilde{P}_{\text{up}}$ | $-\varepsilon + \tilde{P}_{\text{up}}$ | 1 | 1 | 1 | 0 | 0 | 0 | 0 | 0 |
| $u_t^8$ | $2\varepsilon$ | $2\varepsilon + \tilde{P}_{\text{up}}$ | $-\varepsilon + \tilde{P}_{\text{up}}$ | 1 | 1 | 1 | 0 | 0 | 0 | 0 | 0 |
| $u_t^9$ | 0 | 0 | 0 | 0 | 0 | 0 | 0 | 0 | 0 | 0 | 0 |

(2) AFFINELY INDEPENDENT POINTS FOR CONSTRAINTS (59)

| Point | $\tilde{P}_{t-1}$ | $\tilde{P}_t$ | $\tilde{P}_{t+1}$ | $u_{t-1}$ | $u_t$ | $u_{t+1}$ | $s_t$ | $s_{t+1}$ | $d_t$ | $d_{t+1}$ | $\mathcal{T}_t^3$ |
|---|---|---|---|---|---|---|---|---|---|---|---|
| $u_t^1$ | 0 | $\tilde{P}_{\text{start}}$ | $\tilde{P}_{\text{start}} + \tilde{P}_{\text{up}}$ | 0 | 1 | 1 | 1 | 0 | 0 | 0 | 0 |
| $u_t^2$ | 0 | 0 | $\tilde{P}_{\text{start}}$ | 0 | 0 | 1 | 0 | 1 | 0 | 0 | 0 |
| $u_t^3$ | 0 | 0 | 0 | 1 | 0 | 0 | 0 | 0 | 1 | 0 | 0 |
| $u_t^4$ | 0 | $\min(\tilde{P}_{\text{up}}, \tilde{P}_{\text{shut}})$ | 0 | 1 | 1 | 0 | 0 | 0 | 0 | 1 | 0 |
| $u_t^5$ | 0 | $\min(\tilde{P}_{\text{start}}, \tilde{P}_{\text{shut}})$ | 0 | 0 | 1 | 0 | 1 | 0 | 0 | 1 | 1 |
| $u_t^6$ | 0 | $\tilde{P}_{\text{up}}$ | $2\tilde{P}_{\text{up}}$ | 1 | 1 | 1 | 0 | 0 | 0 | 0 | 0 |
| $u_t^7$ | $\varepsilon$ | $\varepsilon + \tilde{P}_{\text{up}}$ | $\varepsilon + 2\tilde{P}_{\text{up}}$ | 1 | 1 | 1 | 0 | 0 | 0 | 0 | 0 |
| $u_t^8$ | 0 | 0 | 0 | 0 | 0 | 0 | 0 | 0 | 0 | 0 | 0 |
| $u_t^9$ | 0 | 0 | $\tilde{P}_{\text{start}}$ | 1 | 0 | 1 | 0 | 1 | 1 | 0 | 0 |

(3) AFFINELY INDEPENDENT POINTS FOR CONSTRAINTS (60)

| Point | $\tilde{P}_{t-1}$ | $\tilde{P}_t$ | $\tilde{P}_{t+1}$ | $u_{t-1}$ | $u_t$ | $u_{t+1}$ | $s_t$ | $s_{t+1}$ | $d_t$ | $d_{t+1}$ | $\mathcal{T}_t^3$ |
|---|---|---|---|---|---|---|---|---|---|---|---|
| $u_t^1$ | 0 | $\tilde{P}_{\text{start}}$ | $\tilde{P}_{\text{start}} + \tilde{P}_{\text{up}}$ | 0 | 1 | 1 | 1 | 0 | 0 | 0 | 0 |
| $u_t^2$ | 0 | 0 | $\tilde{P}_{\text{start}}$ | 0 | 0 | 1 | 0 | 1 | 0 | 0 | 0 |
| $u_t^3$ | 0 | 0 | 0 | 1 | 0 | 0 | 0 | 0 | 1 | 0 | 0 |
| $u_t^4$ | 0 | 0 | 0 | 1 | 1 | 0 | 0 | 0 | 0 | 1 | 0 |
| $u_t^5$ | 0 | 0 | 0 | 0 | 1 | 0 | 1 | 0 | 0 | 1 | 1 |
| $u_t^6$ | 0 | 0 | $\tilde{P}_{\text{up}}$ | 1 | 1 | 1 | 0 | 0 | 0 | 0 | 0 |
| $u_t^7$ | 0 | $\varepsilon$ | $\varepsilon + \tilde{P}_{\text{up}}$ | 1 | 1 | 1 | 0 | 0 | 0 | 0 | 0 |
| $u_t^8$ | $\varepsilon$ | $\varepsilon + \tilde{P}_{\text{up}}$ | $\varepsilon + 2\tilde{P}_{\text{up}}$ | 1 | 1 | 1 | 0 | 0 | 0 | 0 | 0 |
| $u_t^9$ | 0 | 0 | 0 | 0 | 0 | 0 | 0 | 0 | 0 | 0 | 0 |

(4) AFFINELY INDEPENDENT POINTS FOR CONSTRAINTS (61)

| Point | $\tilde{P}_{t-1}$ | $\tilde{P}_t$ | $\tilde{P}_{t+1}$ | $u_{t-1}$ | $u_t$ | $u_{t+1}$ | $s_t$ | $s_{t+1}$ | $d_t$ | $d_{t+1}$ | $\mathcal{T}_t^3$ |
|---|---|---|---|---|---|---|---|---|---|---|---|
| $u_t^1$ | 0 | $\min(\tilde{P}_{\text{start}}, \tilde{P}_{\text{down}})$ | 0 | 0 | 1 | 1 | 1 | 0 | 0 | 0 | 0 |
| $u_t^2$ | 0 | 0 | 0 | 0 | 0 | 1 | 0 | 1 | 0 | 0 | 0 |
| $u_t^3$ | 0 | 0 | 0 | 1 | 0 | 0 | 0 | 0 | 1 | 0 | 0 |
| $u_t^4$ | 0 | $\tilde{P}_{\text{shut}}$ | 0 | 1 | 1 | 0 | 0 | 0 | 0 | 1 | 0 |
| $u_t^5$ | 0 | $\min(\tilde{P}_{\text{start}}, \tilde{P}_{\text{shut}})$ | 0 | 0 | 1 | 0 | 1 | 0 | 0 | 1 | 1 |
| $u_t^6$ | $\tilde{P}_{\text{down}}$ | $\tilde{P}_{\text{down}}$ | 0 | 1 | 1 | 1 | 0 | 0 | 0 | 0 | 0 |
| $u_t^7$ | $\varepsilon + \tilde{P}_{\text{down}}$ | $\varepsilon + \tilde{P}_{\text{down}}$ | $\varepsilon$ | 1 | 1 | 1 | 0 | 0 | 0 | 0 | 0 |
| $u_t^8$ | $\varepsilon + \tilde{P}_{\text{down}}$ | $2\varepsilon + \tilde{P}_{\text{down}}$ | $2\varepsilon$ | 1 | 1 | 1 | 0 | 0 | 0 | 0 | 0 |
| $u_t^9$ | 0 | 0 | 0 | 0 | 0 | 0 | 0 | 0 | 0 | 0 | 0 |

(5) Affinely independent points for constraints (62)

| Point | $\tilde{P}_{t-1}$ | $\tilde{P}_t$ | $\tilde{P}_{t+1}$ | $u_{t-1}$ | $u_t$ | $u_{t+1}$ | $s_t$ | $s_{t+1}$ | $d_t$ | $d_{t+1}$ | $\mathcal{T}_t^3$ |
|---|---|---|---|---|---|---|---|---|---|---|---|
| $u_t^1$ | 0 | 0 | 0 | 0 | 1 | 1 | 1 | 0 | 0 | 0 | 0 |
| $u_t^2$ | 0 | 0 | 0 | 0 | 0 | 1 | 0 | 1 | 0 | 0 | 0 |
| $u_t^3$ | $\tilde{P}_{\text{shut}}$ | 0 | 0 | 1 | 0 | 0 | 0 | 0 | 1 | 0 | 0 |
| $u_t^4$ | $\tilde{P}_{\text{shut}} + \tilde{P}_{\text{down}}$ | $\tilde{P}_{\text{shut}}$ | 0 | 1 | 1 | 0 | 0 | 0 | 0 | 1 | 0 |
| $u_t^5$ | 0 | $\min(\tilde{P}_{\text{start}}, \tilde{P}_{\text{shut}})$ | 0 | 0 | 1 | 0 | 1 | 0 | 0 | 1 | 1 |
| $u_t^6$ | $2\tilde{P}_{\text{down}}$ | $\tilde{P}_{\text{down}}$ | 0 | 1 | 1 | 1 | 0 | 0 | 0 | 0 | 0 |
| $u_t^7$ | $\varepsilon + 2\tilde{P}_{\text{down}}$ | $\varepsilon + \tilde{P}_{\text{down}}$ | $\varepsilon$ | 1 | 1 | 1 | 0 | 0 | 0 | 0 | 0 |
| $u_t^8$ | 0 | 0 | 0 | 0 | 0 | 0 | 0 | 0 | 0 | 0 | 0 |
| $u_t^9$ | $\tilde{P}_{\text{shut}}$ | 0 | 0 | 1 | 0 | 1 | 0 | 1 | 1 | 0 | 0 |

(6) Affinely independent points for constraints (63)

| Point | $\tilde{P}_{t-1}$ | $\tilde{P}_t$ | $\tilde{P}_{t+1}$ | $u_{t-1}$ | $u_t$ | $u_{t+1}$ | $s_t$ | $s_{t+1}$ | $d_t$ | $d_{t+1}$ | $\mathcal{T}_t^3$ |
|---|---|---|---|---|---|---|---|---|---|---|---|
| $u_t^1$ | 0 | 0 | 0 | 0 | 1 | 1 | 1 | 0 | 0 | 0 | 0 |
| $u_t^2$ | 0 | 0 | 0 | 0 | 0 | 1 | 0 | 1 | 0 | 0 | 0 |
| $u_t^3$ | $\tilde{P}_{\text{shut}}$ | 0 | 0 | 1 | 0 | 0 | 0 | 0 | 1 | 0 | 0 |
| $u_t^4$ | $\tilde{P}_{\text{down}}$ | 0 | 0 | 1 | 1 | 0 | 0 | 0 | 0 | 1 | 0 |
| $u_t^5$ | 0 | 0 | 0 | 0 | 1 | 0 | 1 | 0 | 0 | 1 | 1 |
| $u_t^6$ | $\tilde{P}_{\text{down}}$ | 0 | 0 | 1 | 1 | 1 | 0 | 0 | 0 | 0 | 0 |
| $u_t^7$ | $\varepsilon + \tilde{P}_{\text{down}}$ | $\varepsilon$ | 0 | 1 | 1 | 1 | 0 | 0 | 0 | 0 | 0 |
| $u_t^8$ | $2\varepsilon + 2\tilde{P}_{\text{down}}$ | $2\varepsilon$ | $\varepsilon$ | 1 | 1 | 1 | 0 | 0 | 0 | 0 | 0 |
| $u_t^9$ | 0 | 0 | 0 | 0 | 0 | 0 | 0 | 0 | 0 | 0 | 0 |

APPENDIX B.

(50) can be transformed to be the following inequalities:
$$\tilde{P}_{t-1} \leq u_{t-1} + d_t(\tilde{P}_{\text{shut}} - 1) + d_{t+1}(\tilde{P}_{\text{down}} + \tilde{P}_{\text{shut}} - 1), i \in \mathcal{I}^{\geq 2} \quad (73)$$

$$\begin{cases} \tilde{P}_{t-1} \leq u_{t-1} + d_t(\tilde{P}_{\text{shut}} - 1) \quad (52) \\ \tilde{P}_{t-1} \leq u_{t-1} + s_t(1 - \tilde{P}_{\text{down}} - \tilde{P}_{\text{shut}}) + d_t(\tilde{P}_{\text{shut}} - 1) + d_{t+1}(\tilde{P}_{\text{down}} + \tilde{P}_{\text{shut}} - 1) \end{cases}, i \in \mathcal{I}^1 \quad (74)$$

(51) can be transformed to be the following inequalities:
$$\tilde{P}_{t+1} \leq u_{t+1} + s_t(\tilde{P}_{\text{start}} + \tilde{P}_{\text{up}} - 1) + s_{t+1}(\tilde{P}_{\text{start}} - 1), i \in \mathcal{I}^{\geq 2} \quad (75)$$

$$\begin{cases} \tilde{P}_{t+1} \leq u_{t+1} + s_{t+1}(\tilde{P}_{\text{start}} - 1) \quad (53) \\ \tilde{P}_{t+1} \leq u_{t+1} + d_{t+1}(1 - \tilde{P}_{\text{start}} - \tilde{P}_{\text{up}}) + s_t(\tilde{P}_{\text{start}} + \tilde{P}_{\text{up}} - 1) + s_{t+1}(\tilde{P}_{\text{start}} - 1) \end{cases}, i \in \mathcal{I}^1 \quad (76)$$

(59) can be transformed to be the following inequalities:
$$\tilde{P}_{t+1} - \tilde{P}_{t-1} \leq u_{t+1} 2\tilde{P}_{\text{up}} + s_t(\tilde{P}_{\text{start}} - \tilde{P}_{\text{up}}) + s_{t+1}(\tilde{P}_{\text{start}} - 2\tilde{P}_{\text{up}}), i \in \mathcal{I}^{\geq 2} \quad (77)$$

$$\begin{cases} \tilde{P}_{t+1} - \tilde{P}_{t-1} \leq u_{t+1} 2\tilde{P}_{\text{up}} + s_{t+1}(\tilde{P}_{\text{start}} - 2\tilde{P}_{\text{up}}) \\ \tilde{P}_{t+1} - \tilde{P}_{t-1} \leq d_{t+1}(\tilde{P}_{\text{up}} - \tilde{P}_{\text{start}}) + u_{t+1} 2\tilde{P}_{\text{up}} + s_t(\tilde{P}_{\text{start}} - \tilde{P}_{\text{up}}) + s_{t+1}(\tilde{P}_{\text{start}} - 2\tilde{P}_{\text{up}}) \end{cases}, i \in \mathcal{I}^1 \text{ and } \tilde{P}_{\text{up}} - \tilde{P}_{\text{start}} \geq 0 \quad (78)$$

$$\begin{cases} \tilde{P}_{t+1} - \tilde{P}_{t-1} \leq (d_{t+1} - u_t)(\tilde{P}_{\text{up}} - \tilde{P}_{\text{start}}) + u_{t+1} 2\tilde{P}_{\text{up}} + s_{t+1}(\tilde{P}_{\text{start}} - 2\tilde{P}_{\text{up}}) \\ \tilde{P}_{t+1} - \tilde{P}_{t-1} \leq u_{t+1} 2\tilde{P}_{\text{up}} + s_t(\tilde{P}_{\text{start}} - \tilde{P}_{\text{up}}) + s_{t+1}(\tilde{P}_{\text{start}} - 2\tilde{P}_{\text{up}}) \end{cases}, \in \mathcal{I}^1 \text{ and } \tilde{P}_{\text{up}} - \tilde{P}_{\text{start}} < 0 \quad (79)$$

(61) can be transformed to be the following inequalities:
$$\tilde{P}_t - \tilde{P}_{t+1} \leq u_t \tilde{P}_{\text{down}} - s_t[\tilde{P}_{\text{down}} - \tilde{P}_{\text{start}}]^+ + d_{t+1}(\tilde{P}_{\text{shut}} - \tilde{P}_{\text{down}}), i \in \mathcal{I}^{\geq 2} \quad (80)$$

$$\begin{cases} \tilde{P}_t - \tilde{P}_{t+1} \leq -s_t[\tilde{P}_{\text{shut}} - \tilde{P}_{\text{start}}]^+ + u_t \tilde{P}_{\text{down}} + d_{t+1}(\tilde{P}_{\text{shut}} - \tilde{P}_{\text{down}}) \\ \tilde{P}_t - \tilde{P}_{t+1} \leq u_t \tilde{P}_{\text{down}} - s_t[\tilde{P}_{\text{down}} - \tilde{P}_{\text{start}}]^+ + d_{t+1}(\min\{\tilde{P}_{\text{shut}}, \tilde{P}_{\text{start}}\} - \min\{\tilde{P}_{\text{down}}, \tilde{P}_{\text{start}}\}) \end{cases}, i \in \mathcal{I}^1 \text{ and } \varrho_i^1 \geq 0 \quad (81)$$

$$\begin{cases} \tilde{P}_t - \tilde{P}_{t+1} \leq u_t(\min\{\tilde{P}_{\text{down}}, \tilde{P}_{\text{start}}\} + [\tilde{P}_{\text{shut}} - \tilde{P}_{\text{start}}]^+) - s_t[\tilde{P}_{\text{shut}} - \tilde{P}_{\text{start}}]^+ + \\ \qquad d_{t+1}(\min\{\tilde{P}_{\text{shut}}, \tilde{P}_{\text{start}}\} - \min\{\tilde{P}_{\text{down}}, \tilde{P}_{\text{start}}\}) \qquad , i \in \mathcal{I}^1 \text{ and } \varrho_i^1 < 0 \\ \tilde{P}_t - \tilde{P}_{t+1} \leq u_t \tilde{P}_{\text{down}} - s_t[\tilde{P}_{\text{down}} - \tilde{P}_{\text{start}}]^+ + d_{t+1}(\tilde{P}_{\text{shut}} - \tilde{P}_{\text{down}}) \end{cases} \quad (82)$$

where $\varrho_i^1 = [\tilde{P}_{\text{down}} - \tilde{P}_{\text{start}}]^+ - [\tilde{P}_{\text{shut}} - \tilde{P}_{\text{start}}]^+$.

(62) can be transformed to be the following inequalities:
$$\tilde{P}_{t-1} - \tilde{P}_{t+1} \leq u_{t-1} 2\tilde{P}_{\text{down}} - d_t(\tilde{P}_{\text{shut}} - 2\tilde{P}_{\text{down}}) + d_{t+1}(\tilde{P}_{\text{shut}} - \tilde{P}_{\text{down}}), i \in \mathcal{I}^{\geq 2} \quad (83)$$

$$\begin{cases} \tilde{P}_{t-1} - \tilde{P}_{t+1} \leq u_{t-1} 2\tilde{P}_{\text{down}} - d_t(\tilde{P}_{\text{shut}} - 2\tilde{P}_{\text{down}}) + (s_t - d_{t+1})(\tilde{P}_{\text{down}} - \tilde{P}_{\text{shut}}) \\ \tilde{P}_{t-1} - \tilde{P}_{t+1} \leq u_{t-1} 2\tilde{P}_{\text{down}} - d_t(\tilde{P}_{\text{shut}} - 2\tilde{P}_{\text{down}}) \end{cases}, i \in \mathcal{I}^1 \text{ and } \tilde{P}_{\text{down}} - \tilde{P}_{\text{shut}} \geq 0 \quad (84)$$

$$\begin{cases} \tilde{P}_{t-1} - \tilde{P}_{t+1} \leq (s_t - u_t)(\tilde{P}_{\text{down}} - \tilde{P}_{\text{shut}}) + u_{t-1} 2\tilde{P}_{\text{down}} - d_t(\tilde{P}_{\text{shut}} - 2\tilde{P}_{\text{down}}) \\ \tilde{P}_{t-1} - \tilde{P}_{t+1} \leq u_{t-1} 2\tilde{P}_{\text{down}} - d_t(\tilde{P}_{\text{shut}} - 2\tilde{P}_{\text{down}}) + d_{t+1}(\tilde{P}_{\text{shut}} - \tilde{P}_{\text{down}}) \end{cases}, i \in \mathcal{I}^1 \text{ and } \tilde{P}_{\text{down}} - \tilde{P}_{\text{shut}} < 0 \quad (85)$$

(46) can be transformed to be the following inequalities:
$$\tilde{P}_t \leq u_t - s_t(1 - \tilde{P}_{\text{start}}) - d_{t+1}(1 - \tilde{P}_{\text{shut}}), \quad i \in \mathcal{I}^{\geq 2} \quad (86)$$

$$\begin{cases} \tilde{P}_t \leq u_t - s_t[\tilde{P}_{\text{shut}} - \tilde{P}_{\text{start}}]^+ - d_{t+1}(1 - \tilde{P}_{\text{shut}}) \\ \tilde{P}_t \leq u_t - s_t(1 - \tilde{P}_{\text{start}}) - d_{t+1}[\tilde{P}_{\text{start}} - \tilde{P}_{\text{shut}}]^+ \end{cases} \in \mathcal{I}^1 \quad (87)$$

It is obvious that $(87)^1 \gtrsim (74)^1$, $(87)^2 \gtrsim (76)^1$.